\documentclass[a4paper,11pt,reqno]{amsart}
\usepackage[centering,includeheadfoot,margin=2.5 cm]{geometry}
\usepackage{graphics,enumitem,epsfig,textcomp}
\usepackage{amsfonts,amsmath,amssymb,euscript,color,mathrsfs}
\usepackage[utf8]{inputenc}	
\usepackage{hyperref}
\usepackage{soul}

\usepackage[caption=false,subrefformat=parens,labelformat=parens]{subfig}

\usepackage{url}

\usepackage[symbol]{footmisc}
\usepackage{bm}

\numberwithin{equation}{section}
\newtheorem{theorem}{Theorem}[section]
\newtheorem{definition}[theorem]{Definition}
\newtheorem{lemma}[theorem]{Lemma}
\newtheorem{proposition}[theorem]{Proposition}
\newtheorem{corollary}[theorem]{Corollary}
\newtheorem{remark}[theorem]{Remark}
\newtheorem{assumption}[theorem]{Assumption}

\def\eps{\varepsilon}
\def\N{\mathbb{N}}
\def\R{\mathbb{R}}
\def\C{\hbox{\rlap{\kern.24em\raise.1ex\hbox
			{\vrule height1.3ex width.9pt}}C}}
\def\P{\hbox{\rlap{I}\kern.16em P}}
\def\Q{\hbox{\rlap{\kern.24em\raise.1ex\hbox
			{\vrule height1.3ex width.9pt}}Q}}
\def\M{\hbox{\rlap{I}\kern.16em\rlap{I}M}}
\def\Z{\hbox{\rlap{Z}\kern.20em Z}}
\def\({\begin{eqnarray}}
\def\){\end{eqnarray}}
\def\[{\begin{eqnarray*}}
\def\]{\end{eqnarray*}}
\def\part#1#2{\frac{\partial #1}{\partial #2}}

\def\pmb#1{\setbox0=\hbox{$#1$}
	\kern-.025em\copy0\kern-\wd0
	\kern-.05em\copy0\kern-\wd0
	\kern-.025em\raise.0433em\box0 }
\def\bar{\overline}

\def\N{\mathbb{N}}
\def\R{\mathbb{R}}

\def\epsilon{\varepsilon}

\def\E{\mathcal{E}}
\def\P{\mathbb{P}}
\def\Q{\mathbb{Q}}

\def\CL{{\mathrm{CL}}}
\def\TL{{\mathrm{TL}}}
\def\L{{\mathrm{L}}}
\def\W{{\mathrm{W}}}
\def\C{{\mathrm{C}}}

\def\TV{\mathrm{TV}}
\def\Glim{\Gamma\text{-}\lim}

\newcommand*\di{\mathop{}\!\mathrm{d}}
\newcommand\bl{\left(}
\newcommand\br{\right)}

\newcommand\bv{\operatorname{BV}}
\newcommand\tv{\mathrm{TV}}
\newcommand\dive{\operatorname{div}}

\newcommand\diam{\operatorname{diam}}

\def\XXint#1#2#3{{\setbox0=\hbox{$#1{#2#3}{\int}$}
		\vcenter{\hbox{$#2#3$}}\kern-.5\wd0}}

\newcommand\nuv{\lambda_{|v|}}
\newcommand\nuvinv{\lambda_{|1-v|}}
\newcommand\nuvn{\lambda_{|v_n|}}
\newcommand\nuvnk{\lambda_{|v_{n_k}|}}
\newcommand\nuvinvn{\lambda_{|1-v_n|}}

\newcommand\nuvtilde{\lambda_{|\tilde{v}|}}
\newcommand\nuvtildeinv{\lambda_{|1-\tilde{v}|}}

\definecolor{darkblue}{rgb}{0.1,0.1,0.6}

\newcommand\I{{\MakeUppercase{\romannumeral 1}}}
\newcommand\II{{\MakeUppercase{\romannumeral 2}}}
\newcommand\III{{\MakeUppercase{\romannumeral 3}}}
\newcommand\IV{{\MakeUppercase{\romannumeral 4}}}
\newcommand\V{{\MakeUppercase{\romannumeral 5}}}
\newcommand\VI{{\MakeUppercase{\romannumeral 6}}}
\newcommand\dist{\mathrm{dist}}
\newcommand\Id{\mathrm{Id}}
\newcommand\GL{\mathrm{GL}}

\begin{document}


	\centerline{{\Large\textbf{$\Gamma$-Convergence of an Ambrosio-Tortorelli approximation scheme}}}
	\vskip 3mm
	\centerline{{\Large\textbf{for image segmentation}}}
	\vskip 7mm	
	
	\centerline{
 {\large Irene Fonseca}\qquad {\large Lisa Maria Kreusser}\qquad {\large Carola-Bibiane Schönlieb} \qquad {\large Matthew Thorpe}
	}
	\vskip 10mm

	\noindent{\bf Abstract.}
	Given an image $u_0$, the aim of minimising the Mumford-Shah functional is to find a decomposition of the image domain into sub-domains and a piecewise smooth approximation $u$ of $u_0$ such that $u$ varies smoothly within each sub-domain. Since the Mumford-Shah functional is highly non-smooth, regularizations such as the Ambrosio-Tortorelli approximation can be considered which is one of the most computationally efficient approximations of the Mumford-Shah functional for image segmentation. While very impressive numerical results have been achieved in a large range of applications when minimising the functional, no analytical results are currently available for minimizers of the functional in the piecewise smooth setting, and this is the goal of this work.
	Our main result is the $\Gamma$-convergence of the Ambrosio-Tortorelli approximation of the Mumford-Shah functional for piecewise smooth approximations.
	This requires the introduction of an appropriate function space. As a consequence of our $\Gamma$-convergence result, we can infer the convergence of minimizers of the respective functionals.

	\vskip 7mm
	\noindent{\bf AMSC:} 	49J45, 49J55, 62H35, 68U10.

	\noindent{\bf Keywords:} $\Gamma$-convergence, Ambrosio-Tortorelli functional, image segmentation.

\section{Introduction}

Due to their volume and complexity, image and video data are among the largest and fastest growing sources of information, and present some of the biggest challenges for data science. Image segmentation, one of the most fundamental and ubiquitous tasks in image analysis, is the process of partitioning an image into disjoint regions with certain characteristics. Typical examples   include image editing (separating foreground from background, merging multiple images), medical applications (segmenting regions with similar grey-scale values), and biological imaging (detecting cancerous cells, finding cells and nuclei). 

Variational models such as the Mumford-Shah model \cite{MumfordShah} are an important tool for image segmentation. In their model, Mumford and  Shah  formulated an energy minimization problem for computing optimal piecewise smooth approximations of a given image. Particular cases of the minimal partition problem,  its extensions and generalizations are proposed in \cite{ChanV99,ChanVeseActiveContours,ChanVeseMultiphaseLevelSet}.

We consider  the image domain to be represented as $\Omega\subset\R^d$ with $d\geq 1$, where $\Omega$ is an interval for $d=1$ and, for example, a rectangle in the plane for $d=2$. By  $u_0\colon \Omega\to \R^m$ with $m\geq 1$, we  denote a given bounded scalar (grey-scale) or vector-valued (colour) image which should be segmented into two regions. Let $C$ be a closed subset in $\Omega$, made up of a finite set of smooth curves, and the length of curves making up $C$ is denoted by $|C|$. We write $|\cdot|$ for the Euclidean norm.

In the segmentation problem proposed by Mumford and Shah \cite{MumfordShah}, the aim is to find a decomposition of $\Omega$ into sub-domains and an optimal piecewise smooth approximation $u$ of $u_0$ such that $u$ varies smoothly within each sub-domain, and rapidly or discontinuously across the boundaries of the sub-domains. This problem is solved by minimizing the energy functional
\begin{align}\label{eq:mumfordshah}
\E^{MS}(u,C):=\int_{\Omega} |u-u_0|^2 \di x +\mu\int_{\Omega\backslash C} |\nabla u|^2 \di x+\nu |C|,
\end{align}
where $\mu,\nu>0$ are fixed parameters, weighting the different terms in the energy functional. If $(u,C)$ is a minimizer of the above functional, then $u$ is an `optimal' piecewise smooth approximation of the initial, possibly noisy image $u_0$, $C$ can be regarded as approximating the edges of $u_0$, and $u$ is smooth outside of $C$, i.e., in $\Omega\backslash C$. Theoretical results on the existence and regularity of minimizers of \eqref{eq:mumfordshah} are provided by Mumford and Shah \cite{MumfordShah}, Morel and Solimini \cite{morel1988,morel1989,morel1994}, and De Giorgi et al. \cite{DeGiorgi1989}.

For proving existence of minimizers based on the direct method from the calculus of variations, it is necessary to find a topology for which the functional is lower semi-continuous, while ensuring compactness of minimizing sequences. However, the last term in \eqref{eq:mumfordshah} is not lower semi-continuous with respect to any compact topology. This motivates the  formulation of \eqref{eq:mumfordshah} proposed by \cite{DeGiorgi1988} and studied in \cite{dalmaso1992}, where the curve $C$ is replaced by the set $J_u$ of jumps of $u$, leading to the weak formulation of \eqref{eq:mumfordshah}
\begin{align}\label{eq:mumfordshahweak}
\E^{wMS}(u,J_u):=\int_{\Omega} |u-u_0|^2 \di x +\mu\int_{\Omega\backslash J_u} |\nabla u|^2 \di x+\nu |J_u|.
\end{align}
 A constructive existence result for piecewise constant functions $u$ in \eqref{eq:mumfordshahweak} is provided in \cite{morel1988,morel1989}, and a practical multi-scale algorithm based on regions growing and merging is suggested for this case in \cite{koepfler1994}. Ambrosio and Tortorelli proposed two elliptic approximations by $\Gamma$-convergence  \cite{Ambrosio,AmbrosioTortorelli} to the weak formulation \eqref{eq:mumfordshahweak} of the Mumford-Shah functional. 
Approximation \cite{AmbrosioTortorelli} is more commonly used in practise. For $\epsilon>0$ and $(u,v)\in L^2(\Omega)^2$, it is defined as
 \begin{align}\label{eq:mumfordshahweakat}
 \E_\epsilon^{AT}(u,v):=\begin{cases}\int_{\Omega} |u-u_0|^2 \di x +\mu\int_{\Omega} v^2 |\nabla u|^2 \di x+\nu \int_\Omega \bl \epsilon|\nabla v|^2 +\frac{(v-1)^2}{4\epsilon} \br\di x, \\ \hspace*{7.4 cm}  (u,v)\in W^{1,2}(\Omega)^2~\text{with}~0\leq v\leq 1,\\
 +\infty,\hspace*{6.57 cm} \text{otherwise}.
 \end{cases}
 \end{align}
A minimizer $(u,J_u)$ of $\E^{wMS}(u,J_u)$ is approximated by a pair $(u_\epsilon,v_\epsilon)$  of  smooth functions, such that $u_\epsilon\to u$ and $v_\epsilon\to 1$ in the $L^2(\Omega)$-topology as $\epsilon\to 0$ and $v_\epsilon$ is different from 1 only in a small neighbourhood of $J_u$ which shrinks as $\epsilon\to 0$. These elliptic approximations result in a coupled system of two equations with unknowns $u_\epsilon$ and $v_\epsilon$ which can be solved by applying standard numerical methods for PDEs. Further approximations and numerical results are provided in \cite{Chambolle1995,Chambolle1999,march1992}. An approximation by $\Gamma$-convergence to the weak formulation of \eqref{eq:mumfordshah}, based on the finite element method, is discussed in \cite{ChambolleDalMaso1999}. However,  most of the methods for solving the weak formulation of the Mumford-Shah functional \eqref{eq:mumfordshah} do not explicitly compute the partition of the image and the set of curves $C$.

\begin{figure}[htbp]
	\subfloat[Input image]{\includegraphics[width=0.32\textwidth]{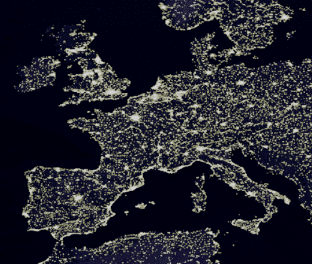}}		
	\subfloat[$\nu=0.2$]{\includegraphics[width=0.32\textwidth]{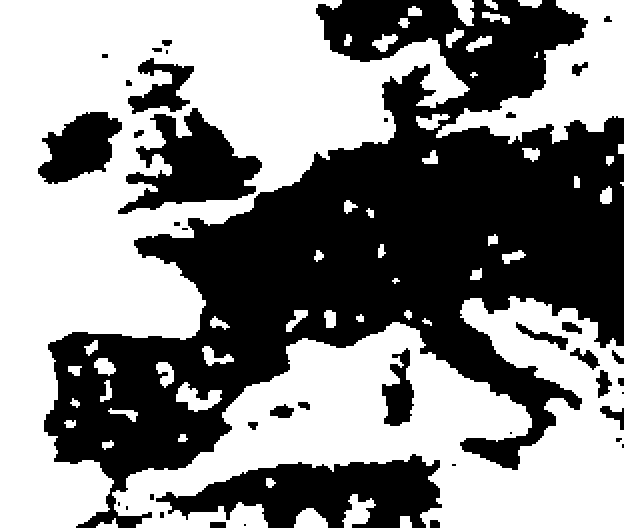}}
	\subfloat[$\nu=0.6$]{\includegraphics[width=0.32\textwidth]{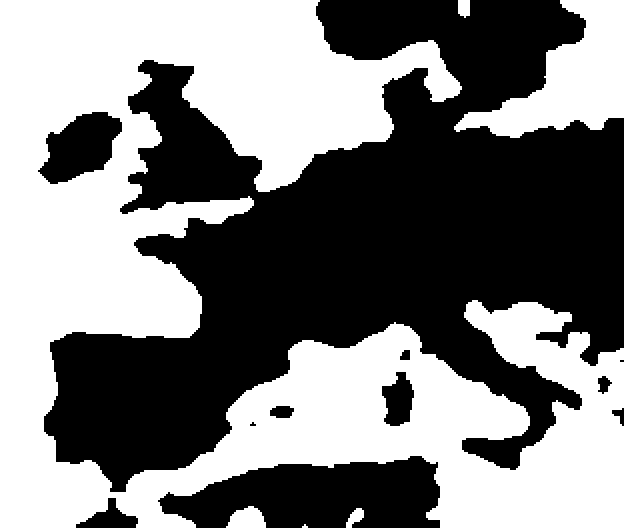}}
	\caption{Image segmentation results for different values for parameter $\nu>0$}\label{fig:segmentation}
\end{figure}

The popular active contour model \cite{ChanVeseActiveContours},
proposed  by Chan and Vese and based on the Mumford-Shah model, can be regarded as a particular case of the Mumford-Shah model \eqref{eq:mumfordshah} by restricting the segmented image $u$ to piecewise constant functions. This model motivates the generalized, widely used multiphase level set model \cite{ChanVeseMultiphaseLevelSet}, also introduced by Chan and Vese. Let $E\subset \Omega$ be an open subset of $\Omega$  inside the boundary curve $C=\partial E$  of length $|C|$, and let $c^{(1)}$ and $c^{(2)}$ be unknown constants. In the active contour model for grey-scale images (i.e., $m=1$), piecewise constant approximations are considered and the energy 
\begin{align*}
\E^{PC}(C,c^{(1)},c^{(2)}):=\int_{E} |c^{(1)}-u_0|^2 \di x +\int_{\Omega\backslash E} (c^{(2)}-u_0)^2 \di x+\nu |C|
\end{align*}
is minimized with respect to $c^{(1)},c^{(2)}\in \R$, and $C$. The parameter  $\nu > 0$ is assumed to be given. The first two terms of $\E^{PC}$ penalize the discrepancy between the input image $u_0$ and its piecewise constant approximation with grey-scale values $c^{(1)}$ in $E$ and $c^{(2)}$ on $\Omega\backslash E$, respectively. The last term controls the regularity of the segmentation by penalizing the length of the boundary curve $C$. Instead of minimizing over all curves $C$, we can represent $C$  implicitly as the zero-crossing of a level set function $\phi\colon \Omega\to \R$, i.e., ${C:=\{x\in\Omega\colon \phi(x)=0\}}$, and we assume that the inside (i.e.\ the set $E$) and the outside (i.e., the set $\Omega\backslash E$) of $C$ are distinguished by positive and negative signs of $\phi$, respectively, to be precise,
\begin{align*}
	\phi(x)>0 \quad \text{in } E,\qquad \phi(x)<0 \quad \text{on } \Omega\backslash E, \qquad \phi(x)=0 \quad \text{on } \partial E.
\end{align*} 
A typical example of a level set function is the signed distance function to the curve.
In its level set formulation, the energy functional can be rewritten as
\begin{align}
\label{eq:PCMumford}
\begin{split}
\E^{lsPC}(\phi,c^{(1)},c^{(2)})&:=\int_{\Omega} |c^{(1)}-u_0|^2 H_\delta(\phi) \di x +\int_{\Omega} |c^{(2)}-u_0|^2 (1-H_\delta(\phi)) \di x \\
& \qquad +\nu\int_{\Omega} |\nabla H_\delta(\phi)|\di x, 
\end{split}
\end{align}
where $H_\delta$ with $\delta>0$ denotes a smooth approximation of the Heaviside function $H$, defined as $H(z)=1$ for $z> 0$ and $H(z)=0$ for $z<0$. Hence, the aim of the active contour model is to find a two-phase segmentation of the image, given by $u(x):=c^{(1)}H_\delta(\phi(x))+c^{(2)}(1-H_\delta(\phi(x)))$, $x\in \Omega$. In Figure~\ref{fig:segmentation}, the segmentation of a given image (based on the implementation in \cite{getreuer}) into two regions, marked in black and white, is shown for $\nu=0.2$ and $\nu=0.6$. The value of the parameter $\nu$ governs the smoothness of the boundary of the segmentation, i.e., for larger values of $\nu$ the interface between white and black areas becomes smaller. This example also illustrates how crucial the parameter choice in this class of models is. 

Following the level set approach, piecewise smooth segmentations are considered in \cite{Tsai2001,ChanVeseMultiphaseLevelSet} by replacing the constants  $c^{(1)},c^{(2)}$ by smooth functions in $E$ and on $\Omega\backslash E$, respectively. The proposed model can be easily extended to vector-valued functions, such as colour images as in \cite{ChanVeseActiveContours}, for instance.  
Based on the Mumford-Shah functional, this 
leads to the energy functional 
\begin{align}\label{eq:energymultiphaseorig}
\begin{split}
\E^{lsPS}(\phi,c^{(1)},c^{(2)})&:=\int_{\Omega} |c^{(1)}-u_0|^2 H_\delta(\phi) \di x +\int_{\Omega} |c^{(2)}-u_0|^2 (1-H_\delta(\phi)) \di x \\&\quad+\mu\int_{\Omega}\left(|\nabla c^{(1)}|^2H_\delta(\phi)+|\nabla c^{(2)}|^2(1-H_\delta(\phi))\right)\di x +\nu\int_{\Omega} |\nabla H_\delta(\phi)|\di x
\end{split}
\end{align}
for piecewise smooth functions $c^{(1)},c^{(2)}$, proposed independently by Vese and Chan \cite{ChanVeseMultiphaseLevelSet}, and Tsai et al.\ \cite{Tsai2001}. Here,  the regularity of $c^{(1)}$ and $c^{(2)}$ is controlled by the  parameter $\mu>0$, and the smoothness of the boundary of the segmentation is governed by $\nu>0$. Numerical results have been obtained independently and contemporaneously by Vese and Chan \cite{ChanVeseMultiphaseLevelSet} and Tsai et al.\ \cite{Tsai2001}. These results show that piecewise smooth regions can be reconstructed very well by the model,  that jumps are well located and without smearing, and that the piecewise constant case can be recovered.

In what follows, we want to study  \eqref{eq:energymultiphaseorig} and its piecewise constant version \eqref{eq:PCMumford}. In particular, the regularity of the piecewise smooth functions $c^{(1)},c^{(2)}$ in \eqref{eq:energymultiphaseorig} is controlled by the parameter $\mu$, and for $\mu\to\infty$ we expect $c^{(1)},c^{(2)}$ to be piecewise constant. This motivates us to study the dependence of the energy  on $\mu$. In addition, it is desirable to control the smoothness of the vector-valued approximations $c^{(1)},c^{(2)}\colon \Omega \to \R^m$, using a parameter $1<p<+\infty$.

The mathematical analysis of  \eqref{eq:energymultiphaseorig}, however, is a highly non-trivial task due to the dependence of the functional on the level set function $\phi$ and on the approximation $H_\delta$ of the non-smooth Heaviside function $H$ as $\phi$ is only implicitly defined and the non-smoothness of $H$ causes difficulties estimating the last term of \eqref{eq:energymultiphaseorig}. They also render the numerical minimization more difficult. To get around this, we propose another formulation that is more amendable to mathematical analysis.
 Since the Heaviside function $H$ only takes values in $\{0,1\}$, this suggests to replace $H(\phi)$ by an indicator function $v$.
These considerations lead to the energy functional  
\begin{align}\label{eq:energymultiphase}
\begin{split}
\E^{PS}_\mu(v,c^{(1)},c^{(2)})&:=\int_{\Omega} |c^{(1)}-u_0|^p |v| \di x +\int_{\Omega} |c^{(2)}-u_0|^p |1-v| \di x \\&\quad+\mu\int_{\Omega} \left(|\nabla c^{(1)}|^p |v|+|\nabla c^{(2)}|^p|1-v|\right)\di x +\nu\int_{\Omega} |\nabla v|\di x,
\end{split}
\end{align}
in place of \eqref{eq:energymultiphaseorig}. For $v=\chi_E$ for some measurable set $E$ with finite perimeter $\operatorname{Per}(E;\Omega)$,  \eqref{eq:energymultiphase} may be equivalently written as
\begin{align*}
    \E^{PS}_\mu(v,c^{(1)},c^{(2)})&=\int_{E}\left( |c^{(1)}-u_0|^p +\mu |\nabla c^{(1)}|^p\right) \di x +\int_{\Omega\backslash E} \left( |c^{(2)}-u_0|^p +
    \mu |\nabla c^{(2)}|^p\right)\di x \\&\quad +\nu \operatorname{Per}(E;\Omega).
\end{align*}

To overcome the non-smoothness of the last term of \eqref{eq:energymultiphase}, several regularization methods and approximations have been proposed in the literature for the numerical minimization. One of the most computationally efficient approximations  of the Mumford-Shah functional was proposed by Ambrosio and Tortorelli \cite{Ambrosio,AmbrosioTortorelli}, and uses the Ginzburg-Landau functional $\E_{\epsilon}^{GL}$ 
defined as
\begin{align}\label{eq:energyginzburglandau}
\E_{\epsilon}^{GL}(v):=\int_{\Omega} \left(\epsilon|\nabla v|^2+\frac{1}{\epsilon}W(v)\right) \di x
\end{align}
which generalizes the approximation in \eqref{eq:mumfordshahweakat}.
Here, $\epsilon>0$ is a positive constant, and the function $W\colon \R\to[0,+\infty)$ is a double well potential with wells at $0$ and $1$, satisfying the following assumption.

\begin{assumption}\label{ass:doublewell}
	Let  $W\colon \R\to[0,+\infty)$ be such that 
	\begin{itemize}
		\item $W$ is continuous,
		\item $W(t)=0$ if and only if $t\in\{0,1\}$, and
		\item there exist $L>0$ and $T>0$ such that
		\begin{align}\label{eq:lingrowth}
		W(t)\geq L|t|\quad \text{for all }t\in\R\text{ with }|t|\geq T.
		\end{align}
	\end{itemize}
\end{assumption}
The most common example for $W$ is $W(x):=x^{2}(x-1)^{2}$. The Ginzburg-Landau functional \eqref{eq:energyginzburglandau} plays an important role due to the work of Modica and Mortola \cite{Modica1977a,Modica1977} who proved that the Ginzburg–Landau functional \eqref{eq:energyginzburglandau}  can be used for approximating the $\tv$ energy, the last term in \eqref{eq:energymultiphase}. In the context of image processing, examples of using the Ginzburg-Landau functional  are given by \cite{Chambolle1995,Chambolle1999}, which relate to previous works by Ambrosio and Tortorelli \cite{ambrosio_fusco_pallara,Ambrosio} on diffuse interface approximation models.

The framework \eqref{eq:energymultiphase} is a very powerful, flexible method that can  segment many types of images, including those that are either difficult or impossible to segment with classical thresholding or gradient-based methods. Using appropriate approximations of the non-smooth terms, this model has been implemented successfully, and very impressive numerical results have been achieved in a large range of applications.  However, no analytical results  are currently available for minimizers of \eqref{eq:energymultiphase} in the piecewise smooth setting, and this is the goal of this work.

\subsection{Contributions}
We will prove
$\Gamma$-convergence of an Ambrosio-Tortorelli approximation of \eqref{eq:energymultiphase}, 
\begin{align}
\bar{\E}_{\mu_\epsilon,\epsilon}(v,c^{(1)},c^{(2)})&:=\int_{\Omega} \left(|c^{(1)}-u_0|^p |v| + |c^{(2)}-u_0|^p |1-v|\right) \di x \label{eq:energyambrosio} \\
 &\quad+\mu_{\epsilon}\int_{\Omega}\left(|\nabla c^{(1)}|^p|v|+|\nabla c^{(2)}|^p|1-v|\right)\di x +\frac{\nu}{c_W} \int_{\Omega}\left( \epsilon|\nabla v|^2+\frac{1}{\epsilon}W(v)\right) \di x, \notag
\end{align}
to the  functional \eqref{eq:energymultiphase}, 
where the positive scaling parameter $\mu_\epsilon$ approximates $\mu>0$, $\nu>0$ is another   scaling parameter, and
\begin{align}\label{eq:cw}
c_W:=2\int_0^1 \sqrt{W(t)}\di t>0.
\end{align}  
In particular, minimizers of \eqref{eq:energyambrosio} will converge to minimizers of \eqref{eq:energymultiphase}, giving new insights into numerical methods for determining minimizers of \eqref{eq:energymultiphase}.

Minimizers of \eqref{eq:energyambrosio} correspond to the segmentation of the vector-valued images $u_0 \colon \Omega \to \R^m$ with $m\geq 1$. 
Since the wells of $W$ are at $0$ and $1$, this suggests that $v$ is an indicator function in the limit $\epsilon \to 0$, and the segmentation, consisting of    smooth approximations $c^{(1)},c^{(2)}\colon \Omega \to \R^m$,   is obtained from $v\colon \Omega \to \R$.

For piecewise constant segmentations of the form $c^{(1)} v+c^{(2)} (1-v)= c^{(1)} \chi_E+c^{(2)} \chi_{\Omega\backslash E}$ for $v=\chi_E$ and constants $c^{(1)},c^{(2)}\in\R^m$, the energy functional \eqref{eq:energyambrosio} reduces to $\bar{\E}_\epsilon\colon L^1(\Omega;\R)\times \R^m \times \R^m$, where
\begin{align*}
\begin{split}
\bar{\E}_{\epsilon}(v,c^{(1)},c^{(2)})&:=\int_{\Omega} \left(|c^{(1)}-u_0|^p |v| + |c^{(2)}-u_0|^p |1-v|\right) \di x+\frac{\nu}{c_W} \int_{\Omega} \left(\epsilon|\nabla v|^2+\frac{1}{\epsilon}W(v)\right) \di x
\end{split}
\end{align*}
for $v\in W^{1,2}(\Omega;\R)$, and $\bar{\E}_{\epsilon}(v,c^{(1)},c^{(2)})=+\infty$ otherwise.
As an illustrative example, we prove  $\Gamma$-convergence of $\bar{\E}_{\epsilon_n}$ to $\bar{\E}\colon L^1(\Omega;\R)\times \R^m \times \R^m$, where
\begin{align*}
\bar{\E}(v,c^{(1)},c^{(2)})&=\int_{\Omega}\left( |c^{(1)}-u_0|^p |v| + |c^{(2)}-u_0|^p|1-v|\right) \di x+\nu  \tv(v)
\end{align*} 
for $v\in\bv(\Omega;\{0,1\})$, and $\bar{\E}(v,c^{(1)},c^{(2)})=+\infty$ otherwise.
Here, $\tv(v)$ denotes the total variation of $v$ in $\Omega$.

For piecewise smooth segmentations of the form $c^{(1)} v+c^{(2)}(1-v)$ where the approximations $c^{(1)},c^{(2)}$ are functions, any $\Gamma$-convergence result requires  $c^{(1)}$ and $c^{(2)}$ to be defined only for $x\in\Omega$ for which $v(x) \neq 0$ and $1-v(x) \neq 0$, respectively, where the sets $\{v=0\}$ and $\{v=1\}$ depend on $v$. Given a function $v \in \L^1((\Omega,\mathcal{L}^d\lfloor_\Omega);\R)$, we want $c^{(1)}$ and $c^{(2)}$, defined on $\Omega$, to be $\nuv$- and $\nuvinv$-measurable, respectively. To achieve this, we introduce the space $\CL^p(\Omega)$  in Section\ref{sec:clp}, motivated by the space $\TL^p(\Omega)$ in \cite{GarciaLimit2016}. Denoting the $d$-dimensional Lebesgue measure by $\mathcal{L}^d$, we say that $(v,c^{(1)},c^{(2)})\in \CL^p(\Omega)$ if $v\in \L^1((\Omega,\mathcal{L}^d\lfloor_\Omega);\R)$, $c^{(1)}\in \L^p((\Omega,\nuv);\R^m), c^{(2)}\in \L^p((\Omega,\nuvinv);\R^m)$, 
where $\nuv$ and $\nuvinv$ are defined by
\begin{align}\label{eq:clp_measure}
\nuv := \begin{cases}\frac{|v|}{\|v\|_{\L^1(\Omega;\R)}}  \mathcal{L}^d \lfloor_\Omega, & \|v\|_{\L^1(\Omega;\R)}\neq 0,\\
0~\mathcal{L}^d\text{-a.e.},& \text{otherwise,}
\end{cases}
\quad \nuvinv :=  \begin{cases} \frac{|1-v|}{\|1-v\|_{\L^1(\Omega;\R)}}  \mathcal{L}^d \lfloor_\Omega,  & \|1-v\|_{\L^1(\Omega;\R)}\neq 0,\\
0~\mathcal{L}^d\text{-a.e.},& \text{otherwise.}
 \end{cases}
\end{align}
We denote the space of distributions on $\Omega$ by $\mathcal D'(\Omega)$, and we consider the space 
\begin{align*}
	\L^{1,p}(\Omega) :=\{ f\in\mathcal D'(\Omega) \enspace \colon \enspace \nabla f \in \L^p(\Omega)\}
\end{align*} 
endowed with the seminorm $\|f\|_{\L^{1,p}}:=\|\nabla f\|_{\L^p}$.
The reformulation of the first term in the second line of the energy functional  \eqref{eq:energyambrosio} with $c^{(1)}\in \L^p((\Omega,\nuv);\R^m), c^{(2)}\in \L^p((\Omega,\nuvinv);\R^m)$  requires the definition  of metric measure Sobolev spaces $\L^{1,p}((\Omega,\nuv);\R^m)$, $\L^{1,p}((\Omega,\nuvinv);\R^m)$, with seminorms $\|\cdot\|_{\L^{1,p}(\nuv)}$ and $\|\cdot\|_{\L^{1,p}(\nuvinv)}$, respectively, which are introduced in Section~\ref{sec:metricmeasurespaces}. The Sobolev space $\W^{1,p}((\Omega,\nuv);\R^m)$ is defined by
\begin{align*}
	\W^{1,p}((\Omega,\nuv);\R^m):=\L^{p}((\Omega,\nuv);\R^m)\cap \L^{1,p}((\Omega,\nuv);\R^m).
\end{align*}
Using the notation of metric measure spaces, we consider a rescaled formulation of the energy functional \eqref{eq:energyambrosio}:
\begin{align}\label{eq:energy}
\begin{split}
\E_{\mu_\epsilon,\epsilon}(v,c^{(1)},c^{(2)})&:=
\| c^{(1)}-u_0\|_{\L^{p}(\nuv;\R^m)}^p+\| c^{(2)}-u_0\|_{\L^{p}(\nuvinv;\R^m)}^p
+\mu_{\epsilon}\| c^{(1)}\|^p_{\L^{1,p}(\nuv)}\\&\quad+\mu_{\epsilon}\| c^{(2)}\|^p_{\L^{1,p}(\nuvinv)} +\frac{\nu}{c_W} \int_{\Omega} \left(\epsilon|\nabla v|^2+\frac{1}{\epsilon}W(v) \right) \di x.
\end{split}
\end{align}
We distinguish between  two cases for the limit of the positive scaling parameter $\mu_\epsilon$, namely  $\mu_\epsilon \to \mu$ with $\mu>0$, and $\mu_\epsilon \to +\infty$ as $\epsilon\to 0$.
For  $+\infty>\mu\geq 0$, we define the limit functional of \eqref{eq:energy} by
\begin{align}\label{eq:energylimit}
\begin{split}
\E_\mu(v,c^{(1)},c^{(2)})&=
\| c^{(1)}-u_0\|_{\L^{p}(\nuv;\R^m)}^p+\| c^{(2)}-u_0\|_{\L^{p}(\nuvinv;\R^m)}^p
\\&\qquad+\mu \| c^{(1)}\|^p_{\L^{1,p}(\nuv)}+\mu \| c^{(2)}\|^p_{\L^{1,p}(\nuvinv)}+\nu  \tv(v)
\end{split}
\end{align} 
for any $v=\chi_E\in\bv(\Omega;\{0,1\})$ with $E=\{x\in\Omega\colon v(x)=1\}$, $c^{(1)}\in \W^{1,p}((\Omega,\nuv);\R^m)$ and $c^{(2)}\in \W^{1,p}((\Omega,\nuvinv);\R^m)$, and $\E_\mu(v,c^{(1)},c^{(2)})=+\infty$ otherwise.
Note that for any Lebesgue measurable set $E\subset \Omega $ such that $\chi_E\in\bv(\Omega;\{0,1\})$, 
the limit functional $\E_\mu$ reduces to
\begin{align*}
\E_\mu(\chi_E,c^{(1)},c^{(2)})&=\frac{1}{|E|}\int_{E} |c^{(1)}-u_0|^p \di x+\frac{1}{|\Omega\backslash E|}\int_{\Omega\backslash E} |c^{(2)}-u_0|^p \di x\\&\qquad+\frac{\mu}{|E|} \| c^{(1)}\|^p_{\L^{1,p}(\lambda_{\chi_E})}+\frac{\mu}{|\Omega\backslash E|} \| c^{(2)}\|^p_{\L^{1,p}(\lambda_{\chi_{\Omega\backslash E}}))}  +\nu  \tv(\chi_E), 
\end{align*} 
where $|E|=\mathcal{L}^d(E)$ denotes the $d$-dimensional Lebesgue measure of $E$, and $\lambda_{\chi_E},\lambda_{\chi_{\Omega\backslash E}}$ are defined as in \eqref{eq:clp_measure}. For a bounded domain $E$ with smooth boundary, the norms $\| c^{(1)}\|_{\L^{1,p}(\lambda_{\chi_E})}$ and $\| c^{(1)}\|_{\L^{1,p}(E)}=\| \nabla c^{(1)}\|_{\L^{p}(E)}$ are equivalent.

Our main result is  the $\Gamma$-convergence of the variational model \eqref{eq:energy} to \eqref{eq:energylimit} as $\epsilon\to 0$. 
\begin{theorem}\label{th:mainresult}
Let $\Omega \subset \R^d$ be an open, bounded  set, let $1<p<+\infty$, and let $\E_{\mu_\epsilon,\epsilon}\colon \CL^p(\Omega)\to [0,+\infty]$ be defined by
		\begin{align}\label{eq:approxenergy}
\E_{\mu_\epsilon,\epsilon}(v,c^{(1)},c^{(2)}):=\begin{cases} 
\| c^{(1)}-u_0\|_{\L^{p}(\nuv;\R^m)}^p+\| c^{(2)}-u_0\|_{\L^{p}(\nuvinv;\R^m)}^p
+\mu_{\epsilon}\| c^{(1)}\|^p_{\L^{1,p}(\nuv)}\\\quad+\mu_{\epsilon}\| c^{(2)}\|^p_{\L^{1,p}(\nuvinv)} +\frac{\nu}{c_W} \int_{\Omega}\epsilon|\nabla v|^2+\frac{1}{\epsilon}W(v) \di x,\\ \hspace*{1.7 cm} \text{if } v\in \W^{1,2}((\Omega,\mathcal{L}^d\lfloor_\Omega);\R),c^{(1)}\in \W^{1,p}((\Omega,\nuv);\R^m), \\ \hspace*{1.7 cm}c^{(2)}\in \W^{1,p}((\Omega,\nuvinv);\R^m),\\+\infty,\hspace*{0.85 cm}\text{otherwise.}\end{cases}
\end{align}
Then, the functionals $\E_{\mu_\epsilon,\epsilon}$ $\Gamma$-converge, with respect to the $\CL^p(\Omega)$ topology, 
\begin{align}\label{eq:limitenergycase1}
\E_\mu(v,c^{(1)},c^{(2)}):=\begin{cases} \int_{\Omega} |c^{(1)}-u_0|^p \di \nuv(x)+\int_\Omega |c^{(2)}-u_0|^p \di \nuvinv (x)+\mu \| c^{(1)}\|^p_{\L^{1,p}(\nuv)}\\\quad+\mu \| c^{(2)}\|^p_{\L^{1,p}(\nuvinv)}+\nu  \tv(v),
\\\hspace*{1.7cm} \text{if } v=\chi_E\in\bv(\Omega;\{0,1\}) \text{ for }E:=\{x\in\Omega\colon v(x)=1\},
\\ \hspace*{1.7 cm} c^{(1)}\in \W^{1,p}((\Omega,\nuv);\R^m), c^{(2)}\in \W^{1,p}((\Omega,\nuvinv);\R^m),
\\+\infty,\hspace*{0.85cm} \text{otherwise,}\end{cases}
\end{align}
if  $\mu_\epsilon \to \mu$ with $\mu>0$, as $\epsilon\to 0$, and
\begin{align}\label{eq:limitenergycase2}
\E_\infty(v,c^{(1)},c^{(2)}):=\begin{cases} \int_{\Omega} |c^{(1)}-u_0|^p \di \nuv(x)+\int_\Omega |c^{(2)}-u_0|^p \di \nuvinv (x)+\nu  \tv(v),
\\\hspace*{1.7cm}\text{if }v=\chi_E\in\bv(\Omega;\{0,1\}) \text{ for }E:=\{x\in\Omega\colon v(x)=1\},
\\ \hspace*{1.7 cm} c^{(1)}=c_1~\mathcal{L}^d\text{-a.e. } x\in E, c^{(2)}=c_2~ \mathcal{L}^d\text{-a.e. }  x\in \Omega\backslash E \text{ for  }
\\\hspace*{1.7 cm} \text{constants } c_1,c_2\in\R^m,
\\+\infty,\hspace*{.85cm} \text{otherwise,}\end{cases}
\end{align}
if $\mu_\epsilon \to +\infty$ as $\epsilon\to 0$.
\end{theorem}

Provided that the compactness property holds, i.e.\ every bounded sequence $(v_n,c_n^{(1)},c_n^{(2)})\in \CL^p(\Omega)$ satisfying $\sup_{n\in\N} \E_{\mu_{\epsilon_n},\epsilon_n}(v_n,c^{(1)}_n,c^{(2)}_n)<\infty$ is relatively compact, the convergence of minimizers follows from the $\Gamma$-convergence of the energy functional $\E_{\mu_\epsilon,\epsilon}$. We prove the compactness property in Theorem \ref{th:compactnesssmoothcase1} if there exists $v=\chi_E$ for some Lebesgue measurable set $E\subset \Omega$, if there exist $\kappa>0$, $r_0>0$  such that $P(E;B_r(x)) \geq \kappa r^d$ for every $x \in \partial^* E$,  if there exist ${\kappa>0}$, $r_0>0$  such that $P(\Omega \setminus E;B_r(x)) \geq \kappa r^d$ for every $x \in \partial^* (\Omega \setminus E)$, and if there exists a subsequence $\{v_{n_k}\}$ of $\{v_n\}$  such that ${v_{n_k}\to v}$ in $\L^1(\Omega;\R)$.
In particular, we prove the following corollary.

	\begin{corollary}[Convergence of minimizers]\label{cor:convergenceminimizer}
	Let $\Omega\subset \R^d$ be an open, bounded set  with $d\geq 2$. Suppose that $(v_n,c_n^{(1)},c_n^{(2)})\in \CL^p(\Omega)$ is a minimizer of the energy $\E_{\mu_{\epsilon_n},\epsilon_n}$ in \eqref{eq:approxenergy}, for positive sequences $\{\epsilon_n\},  \{\mu_{\epsilon_n}\}$ with $\lim_{n\to\infty}\epsilon_n= 0$ and $\lim_{n\to \infty} \mu_{\epsilon_n}=\mu\in(0,+\infty]$. If there exists $v=\chi_E$ for some Lebesgue measurable set $E\subset \Omega$, if there exist $\kappa>0$, $r_0>0$  such that $P(E;B_r(x)) \geq \kappa r^d$ for every $x \in \partial^* E$,  if there exist ${\kappa>0}$, $r_0>0$  such that $P(\Omega \setminus E;B_r(x)) \geq \kappa r^d$ for every $x \in \partial^* (\Omega \setminus E)$, and if there exists a subsequence $\{v_{n_k}\}$ of $\{v_n\}$  such that ${v_{n_k}\to v}$ in $\L^1(\Omega;\R)$, then there exists $(v,c^{(1)},c^{(2)})\in \CL^p(\Omega)$ such that, up to a subsequence (not relabeled), $(v_n,c_n^{(1)},c_n^{(2)})$ converges  to $(v,c^{(1)},c^{(2)})$ in $\CL^p(\Omega)$, and $(v,c^{(1)},c^{(2)})$ minimizes the energy $\E_\mu$ in \eqref{eq:limitenergycase1} and \eqref{eq:limitenergycase2} for $\mu<+\infty$ and $\mu=+\infty$, respectively, over $\CL^{p}(\Omega)$.	
\end{corollary} 

While we focus on image segmentations into two segments in this work, the analysis can be extended to images which are partitioned into more than two segments.

\subsection{Overview}
In Section \ref{sec:def}, we give some preliminary material which includes the definition of metric measure spaces, transportation theory, $\Gamma$-convergence and the space $\CL^p$. Section \ref{sec:pw_constant} is devoted to the proof of Theorem \ref{th:mainresult} for piecewise constant segmentations, i.e., $c^{(1)},c^{(2)}\in\R^m$. In Section \ref{sec:pw_smooth}, we prove Theorem \ref{th:mainresult} for piecewise smooth approximations and we show the convergence of minimizers of the respective functionals.

\section{Definitions and preliminary results}\label{sec:def}

\subsection{Notation}
	Throughout this paper, let $\chi_E$ denote the characteristic function of a set $E\subset \R^d$. We write $\mathcal{L}^d$ for the  $d$-dimensional Lebesgue measure on $\R^d$,  and $|E|=\mathcal{L}^d(E)$ stands for the $d$-dimensional Lebesgue measure of $E$. For an open set $\Omega\subset \R^d$, we designate by $\mathcal{B}(\Omega)$ the Borel $\sigma$-algebra on $\Omega$, and by $\mathcal{P}(\Omega)$ the set of Borel probability measures on $\Omega$. For the measure space $(\Omega, \mathcal{B}(\Omega),\lambda)$, where $\lambda$ is a measure on $(\Omega, \mathcal{B}(\Omega))$, we often write $(\Omega,\lambda)$. For the $L^p$ space of all measurable functions from $(\Omega,\lambda)$ to $\R^m$, we write $\L^{p}((\Omega,\lambda);\R^m)$. If the considered spaces or measures are clear, we may use $\L^p(\Omega)$ or $\L^p(\lambda)$ for ease of notation.
	 The space of functions of bounded variation, $\bv(\Omega;\R)$, is defined as the space of all functions $v\in \L^1(\Omega;\R)$ whose distributional first-order partial derivatives are finite signed Radon measures, defined on the Borel $\sigma$-algebra $\mathcal{B}(\Omega;\R)$, i.e., for all $i=1,\ldots, d$, there exists a finite signed measure $v_i\colon \mathcal{B}(\Omega;\R)\to\R$ such that
		\begin{align*}
		\int_{\Omega}v\frac{\partial\Phi}{\partial x_i}\di x=-\int_{\Omega} \Phi\di v_i
		\end{align*}
	for all $\Phi\in C_c^{\infty}(\Omega;\R)$. The measure $v_i$ is called the weak partial derivative of $v$ with respect to $x_i$, and is denoted by $D_i v$. 
	For $v\in \bv(\Omega;\R)$ we set $Dv:=(D_1 v,\ldots,D_d v)$.
		The total variation of $v$ in $\Omega$ for $v\in \L^1_{loc}(\Omega;\R)$ is defined by
		\begin{align*}
		\TV(v):=\sup\left\{\int_{\Omega} v \dive \Phi \di x ~\colon \Phi\in \C_c^{\infty}(\Omega;\R^d),~\|\Phi\|_{\L^{\infty}(\Omega; \R^d)}\leq 1\right\}.
		\end{align*}

	\subsection{Definition of metric measure spaces}\label{sec:metricmeasurespaces}
	
	Sobolev spaces can be defined on metric measure spaces \cite{SobolevSpacesArbitrary,hajlasz2003sobolev,hajlasz2000sobolev}. For completeness, we recall the standard definitions of Sobolev spaces 
	\begin{align*}
	\W^{1,p}(\Omega)&=\{f\in\mathcal D'(\Omega)  \enspace \colon \enspace f\in \L^p(\Omega), \nabla f \in \L^p (\Omega)\},\\
	\L^{1,p}(\Omega) &=\{ f\in\mathcal D'(\Omega)  \enspace \colon \enspace \nabla f \in \L^p(\Omega)\},
	\end{align*}
	where $\Omega\subset \R^d$ is an open set, $1\leq p\leq +\infty$, and $\mathcal{D}'(\Omega)$ denotes the space of distributions on $\Omega$. The space $\W^{1,p}(\Omega)$ is a Banach space when endowed with the norm $\|f\|_{\W^{1,p}}:=\|f\|_{\L^p}+\|\nabla f\|_{\L^p}$, $\L^{1,p}(\Omega)$ is endowed with the seminorm $\|f\|_{\L^{1,p}}:=\|\nabla f\|_{\L^p}$. Note that $\W^{1,p}(\Omega)\neq \L^{1,p}(\Omega)$ in general. 
	
	The definition of  Sobolev spaces strongly relies on the Euclidean structure of the underlying domain $\Omega$. 
	In order to define Sobolev spaces on metric measure spaces, we need to consider a different approach that does not involve derivatives. From \cite[Theorem 2.2]{hajlasz2003sobolev}, we obtain:
	\begin{theorem}\label{th:relationdiffmetricmeasure}
		Let  $\Omega\subset \R^d$ be a bounded domain with smooth boundary, and let $1<p<+\infty$. Then $f\in \W^{1,p}(\Omega)$,  if and only if $f \in \L^p(\Omega)$, and there is $0 \leq g \in \L^p(\Omega)$  such that
		\begin{align}\label{eq:sobolevproperty}
		|f(x) - f(y)| \leq |x - y|(g(x) + g(y))  \enspace \mathcal{L}^d\text{-a.e.}
		\end{align}
		Moreover, $\|f\|_{\L^{1,p}}$ is equivalent to  $\inf_g \|g\|_{\L^p}$, i.e., there exists a constant $C\geq 1$ such that $\frac{1}{C}\|f\|_{\L^{1,p}} \leq  \inf_g \|g\|_{\L^p}\leq C \|f\|_{\L^{1,p}}$, where the infimum is taken over the class of all
		functions $g$ satisfying~\eqref{eq:sobolevproperty}.
	\end{theorem}

	This definition can be extended to the case in which $\Omega$ is replaced by a metric  space $(\Omega,d)$ equipped with a Borel measure $\lambda$:
	\begin{definition}
		Let $(\Omega,d)$ be a metric space  with a finite positive Borel measure $\lambda$ and finite diameter, $$\diam \Omega:=\sup_{x,y\in \Omega} d(x,y)<+\infty.$$ Let $1< p< +\infty$. The Sobolev spaces $\L^{1,p}(\Omega,d,\lambda)$ and $\W^{1,p}(\Omega,d,\lambda)$ are defined, respectively, as
		\begin{align*}
		\L^{1,p}(\Omega,d,\lambda) :=\{ f\colon \Omega\to \R\enspace \colon\enspace &f \text{ is measurable, and there exist } E\subset \Omega \text{ with } \lambda(E)=0 \text{ and }\\&0\leq g\in \L^p(\lambda) \text{ such that } |f(x)-f(y)|\leq d(x,y)(g(x)+g(y)) \\ &\text{for all }x,y\in \Omega\backslash E \},
		\end{align*}
		and
		\begin{align*}
		\W^{1,p}(\Omega,d,\lambda):=\L^p(\lambda)\cap \L^{1,p}(\Omega,d,\lambda).
		\end{align*}
	\end{definition}
	
	The space $\L^{1,p} (\Omega,d,\lambda)$ is equipped with the seminorm $\|f\|_{\L^{1,p}(\lambda)}:=\inf_g \|g\|_{\L^p(\lambda)}$ where $0\leq g$ satisfies 
		\begin{align}\label{eq:sobolevpropertygenerald}
	|f(x) - f(y)| \leq d(x,y)(g(x) + g(y))  \enspace \lambda\text{-a.e.}
	\end{align}
	 The space $\W^{1,p}(\Omega,d,\lambda)$ is equipped with the norm $\|f\|_{\W^{1,p}(\lambda)}:=\|f\|_{\L^p(\lambda)}+\|f\|_{\L^{1,p}(\lambda)}$. If the metric $d$ is clear, we also write $\L^{1,p} (\Omega,\lambda)$ and $\W^{1,p}(\Omega,\lambda)$. 
	
	\begin{remark}
		Note that other modifications of  $\L^p$ spaces  exist, such as the weighted $\L^p$ space with a weight function $w$ on $\Omega$. However, while these spaces are defined on a domain $\Omega$, we are interested in $\L^p$ spaces, and more generally Sobolev spaces, on some measure space $(\Omega,\lambda)$ for some nonnegative measure $\lambda$. For $f\in \L^{1,p}(\Omega,\lambda)$ where $\lambda=\chi_E \mathcal{L}^d\lfloor_{\Omega}$ is the indicator function of some  measurable bounded domain $E\subset \Omega$  with smooth boundary, we have $f\in \L^{1,p}(E)$. 
		In particular, the norms $\| f\|_{\L^{1,p}(\lambda_{\chi_E})}$ and $\| f\|_{\L^{1,p}(E)}=\| \nabla f\|_{\L^{p}(E)}$ are equivalent.
	\end{remark}

	\subsection{Transportation theory}\label{sec:transportation}
	\begin{definition}
		Let $\Omega\subset \R^d$ be an open set, and let $\lambda,\tilde{\lambda}$ be probability measures on $\Omega$. We define the set of couplings $\Pi(\lambda,\tilde{\lambda})$ between $\lambda$ and $\tilde{\lambda}$ as
		\begin{align*}
		\Pi(\lambda,\tilde{\lambda}):=\{\pi\in \mathcal{P}(\Omega\times \Omega) \colon \pi(E \times \Omega)=\lambda(E),\; \pi(\Omega\times E)=\tilde{\lambda}(E)\; \text{for all measurable}\; E\subset \Omega \}.
		\end{align*}
		The elements $\pi\in\Pi(\lambda,\tilde{\lambda})$ are also referred to as transportation plans between $\lambda$ and $\tilde{\lambda}$.
	\end{definition}	
	\begin{definition}
		Let $1\leq p<+\infty$, $\lambda\in\mathcal{P}(\Omega)$ and $\{\lambda_n\}\subset \mathcal{P}(\Omega)$. A sequence of transportation plans $\{\pi_n\} \subset\Pi(\lambda,\lambda_n)$ is called stagnating if 
		\begin{align}\label{eq:stagnating}
		\lim_{n\to \infty} \int_{\Omega\times \Omega} |x-y|^p \di \pi_n(x,y)=0
		\end{align}
		is satisfied.
	\end{definition}
	Since $\Omega$ is bounded, the existence of a stagnating sequence of transportation plans is equivalent to the weak convergence of probability measures, i.e., $\{\lambda_n\}$ converges weakly-$\ast$ to $\lambda$ if and only if for any $1\leq  p < +\infty$ there is a sequence of transportation plans $\{\pi_n\} \subset\Pi(\lambda,\lambda_n)$ for which \eqref{eq:stagnating} is satisfied \cite{ambrosioGradientFlowsMetric2008,villani}.
	\begin{lemma}\cite{GarciaLimit2016}\label{lem:stagseq_1}
	Let $1\leq p<+\infty$, $\lambda\in\mathcal{P}(\Omega)$,  $\{\lambda_n\}\subset\mathcal{P}(\Omega)$, and let 
	$\{\pi_n\}\subset\Pi(\lambda,\lambda_n)$
	for all $n\in\N$. If $\{\pi_n\}$ is a stagnating sequence  of transportation plans, then
	for any $c \in \L^p((\Omega,\lambda);\R^m)$
	\begin{align*}
		\lim_{n\to\infty}  \int_{\Omega\times \Omega} |c(x)-c(y)|^p \di \pi_n(x,y)=0.
	\end{align*}
\end{lemma}		
\begin{lemma}\cite{GarciaLimit2016}\label{lem:stagseq_2}
	Suppose that the sequence $\{\lambda_n\}$ in $\mathcal{P}(\Omega)$ converges weakly-$\ast$ to $\lambda \in\mathcal{P}(\Omega)$. Let $c_n \in \L^p((\Omega,\lambda_n);\R^m)$, $n\in\N$, and let $c \in \L^p((\Omega,\lambda);\R^m)$.
	Consider two sequences of stagnating transportation plans $\{\pi_n\}$ and $\{\tilde{\pi}_n\}$, with $\pi_n,\tilde{\pi}_n\in \Pi(\lambda,\lambda_n)$. Then,
	\begin{align*}
		\lim_{n\to\infty} \int_{\Omega\times \Omega} |c(x)-c_n(y)|^p \di \pi_n(x,y)=0\iff 	\lim_{n\to\infty}  \int_{\Omega\times \Omega} |c(x)-c_n(y)|^p \di \tilde{\pi}_n(x,y)=0.
	\end{align*}
	\end{lemma}
	\begin{definition}
		Given a Borel map $T\colon \Omega \to \Omega$ and $\lambda \in \mathcal{P}(\Omega)$, the push-forward of $\lambda$ by $T$ is denoted by $ T_\#\lambda\in\mathcal{P}(\Omega)$, and is given by
		\begin{align*}
		T_\# \lambda(E) := \lambda(T^{-1}(E)),\;  E\in \mathcal{B}(\Omega).
		\end{align*}
		\end{definition}
		For any bounded Borel function $\phi\colon \Omega \to \R$, the following change of variables holds:
		\begin{align}\label{eq:changeofvariablegeneral}
		\int_\Omega \phi(x) \di (T_\# \lambda)(x)=\int_\Omega \phi(T(x))\di \lambda(x).
		\end{align}
		\begin{definition}
		A Borel map $T\colon \Omega \to \Omega$ is called a transportation map between the measures $\lambda\in \mathcal{P}(\Omega)$ and $\tilde{\lambda}\in\mathcal{P}(\Omega)$ if $\tilde{\lambda}=T_\# \lambda$.
		\end{definition}
		For a transportation map $T$ between measures $\lambda,\tilde{\lambda}\in\mathcal P(\Omega)$, we associate to $T$ the transportation plan $\pi_T\in\Pi(\lambda,\tilde{\lambda})$ given by
		\begin{align}\label{eq:transportationplan}
		\pi_T:= (\text{Id}\times T)_\# \lambda,
		\end{align}
		where $\text{Id}\times T\colon \Omega \to \Omega \times \Omega$ with $(\text{Id}\times T)(x)=(x,T(x))$. 
		For any $\phi \in \L^1(\Omega\times \Omega,\R)$,
		a change of variables yields
		\begin{align}\label{eq:changeofvariable}
			\int_{\Omega\times \Omega} \phi(x,y) \di \pi_T(x,y)=\int_\Omega \phi(x,T(x))\di \lambda(x).
		\end{align}		
		
		\subsection{\texorpdfstring{$\Gamma$}{Gamma}-convergence}
		We recall the notion of $\Gamma$-convergence  \cite{braides2002,maso1993}.
		\begin{definition}
			Let $(X,d)$ be a metric space, and let $\{\E_n\}$ be a sequence of functions $\E_n\colon X\to[-\infty,+\infty]$. We say that $\{\E_n\}$ $\Gamma$-converges to a function $\E\colon X\to[-\infty,+\infty]$ if the following two properties are satisfied:
			\begin{itemize}
				\item (Liminf inequality) For every $x\in X$ and every sequence $\{x_n\}\subset X$ such that $x_n\to x$ with respect to $d$,
				\begin{align*}
				\E(x)\leq \liminf_{n \to \infty} \E_n(x_n).
				\end{align*}
				\item (Limsup inequality) For every $x\in X$, there exists a sequence $\{x_n\}\subset X$ such that $x_n\to x$  with respect to $d$, and
				\begin{align*}
				\limsup_{n \to \infty} \E_n(x_n)\leq \E(x).
				\end{align*}
			\end{itemize}
			The limit function $\E$ is called the $\Gamma$-limit of the sequence $\{\E_n\}$, and we  write $$\Glim_{n\to\infty} \E_n =\E.$$
		\end{definition}

		\begin{definition}
			Let $(X,d)$ be a metric space. A sequence of nonnegative functionals $\{\E_n\}$ with $\E_n\colon X\to[-\infty,+\infty]$ satisfies the compactness property if for any increasing subsequence $\{n_k\}$ of natural numbers and any bounded sequence $\{x_k\}\subset X$ such that 
		\begin{align*}
			\sup_{k\in\N} \E_{n_k}(x_k)<\infty,
			\end{align*}
			the sequence $\{x_k\}$ is relatively compact in $X$.
			\end{definition}
		
			For functionals $\{\E_n\}$ satisfying the compactness property, the notion of $\Gamma$-convergence is particular useful since it guarantees the convergence of minimizers (or approximations of minimizers) of $\E_n$ to minimizers of $\E$. It also guarantees the convergence of the minimum energy of $\E_n$ to the minimum energy of $\E$. To be precise,
		\begin{proposition}
			Let $\E_n\colon X\to [0,\infty]$ be  nonnegative functionals  not identically equal to $+\infty$, satisfying the compactness property, and $\Gamma$-converging to the functional $\E\colon X\to [0,\infty]$ that is not identically equal to $+\infty$. Then, 
			\begin{align*}
				\lim_{n\to\infty} \inf_{x\in X} \E_n(x)=\min_{x\in X} \E(x).
			\end{align*} 
			Furthermore, every bounded sequence $\{x_n\}_{n\in\N}$ in $X$ for which
			\begin{align}\label{eq:mingammaconv}
				\lim_{n\to\infty}\bl \E_n(x_n)-\inf_{x\in X} \E_n(x) \br=0
			\end{align}
			is relatively compact, and each of its cluster points is a minimizer of $\E$. In particular, if $\E$ has a unique minimizer, then a sequence $\{x_n\}$ satisfying \eqref{eq:mingammaconv} converges to the unique minimizer of~$\E$.
		\end{proposition}

		\begin{theorem}[$\Gamma$-Convergence and Compactness of the Ginzburg-Landau Functional \cite{Modica1977a,Modica1977,Modica1987,Sternberg1988}]\label{th:compactness}
			Let $\Omega\subset\R^d$ be an open, bounded set and $\epsilon_n\to 0$. Suppose that Assumption \ref{ass:doublewell} is satisfied, and define $\E^{GL}_{\epsilon_n}$ by \eqref{eq:energyginzburglandau}.
			Then, $\Glim_{n\to\infty} \E^{GL}_{\epsilon_n} = c_W\TV$, with $c_W$ as in \eqref{eq:cw}.
			Furthermore, let $\{v_n\}\subset \W^{1,2}(\Omega;\R)$ be such that
			\begin{align*}
			M:=\sup_{n\in\N}\E^{GL}_{\epsilon_n}(v_n)<+\infty.
			\end{align*}
			Then there exist a subsequence $\{v_{n_k}\}$ of  $\{v_n\}$ and $v\in \bv(\Omega;\{0,1\})$ such that 
			\begin{align*}
			v_{n_k}\to v\quad \text{in }\L^1(\Omega;\R).
			\end{align*}
		\end{theorem}
		Using the results of Modica and Mortola \cite{Modica1977a,Modica1977}, Modica \cite{Modica1987} and Sternberg \cite{Sternberg1988} independently proved Theorem \ref{th:compactness} under the stronger assumption that 
		\begin{align*}
		\frac{1}{c}|t|^q\leq W(t)\leq c|t|^q
		\end{align*}
		for all $|t|\geq T$ for some $T>0$, $c>0$ and $q\geq 2$.  Fonseca and Tartar \cite{Fonseca} showed that the weaker assumption of linear growth in \eqref{eq:lingrowth} is sufficient for Theorem \ref{th:compactness}.

		\subsection{The space \texorpdfstring{$\CL^p$}{CLp}}\label{sec:clp}
		Let $\Omega\subset \R^d$ be an open set. We define
				\begin{align*}
		\CL^p(\Omega) & := \Bigg\{ (v,c^{(1)},c^{(2)}) \colon v\in \L^1((\Omega,\mathcal{L}^d\lfloor_\Omega);\R), \\
		& \qquad \qquad c^{(1)}\in \L^p((\Omega,\nuv);\R^m), c^{(2)}\in \L^p((\Omega,\nuvinv);\R^m) \Bigg\}, 
		\end{align*}
		where $\nuv$ and $\nuvinv$ are given by \eqref{eq:clp_measure}, i.e.\ $\nuv$ and $\nuvinv$  are probability measures on $\Omega$ which have Lebesgue densities $\frac{|v|}{\|v\|_{\L^1(\Omega;\R)}}$ and $\frac{|1-v|}{\|1-v\|_{\L^1(\Omega;\R)}}$ if $\|v\|_{\L^1(\Omega;\R)}\neq 0$ and $\|1-v\|_{\L^1(\Omega;\R)}\neq 0$, respectively. For $(v,c^{(1)},c^{(2)})$ and $(\tilde{v},\tilde{c}^{(1)},\tilde{c}^{(2)})$ in $\CL^p(\Omega)$, we define the equivalence relation on $\CL^p$ as:
		\begin{align*}
		&(v,c^{(1)},c^{(2)})\sim (\tilde{v},\tilde{c}^{(1)},\tilde{c}^{(2)})\\&\iff  \begin{cases}
		v=\tilde{v}=0 \enspace \mathcal{L}^d\text{-a.e.},\enspace c^{(1)}=\tilde{c}^{(1)} \enspace \nuv\text{-a.e.}, \enspace c^{(2)}=\tilde{c}^{(2)}\enspace\nuvinv\text{-a.e.},
		\\\hspace*{1.7 cm} \text{if }\|v\|_{\L^1}=0 \text{ or } \|\tilde{v}\|_{\L^1}=0,
		\\
		v=\tilde{v}=1 \enspace \mathcal{L}^d\text{-a.e.},\enspace c^{(1)}=\tilde{c}^{(1)} \enspace \nuv\text{-a.e.}, \enspace c^{(2)}=\tilde{c}^{(2)}\enspace\nuvinv\text{-a.e.},
		\\\hspace*{1.7 cm} \text{if }\|1-v\|_{\L^1}=0 \text{ or } \|1-\tilde{v}\|_{\L^1}=0,
		\\
		\nuv=\nuvtilde\enspace \mathcal{L}^d\text{-a.e.},\enspace \nuvinv=\nuvtildeinv \enspace\mathcal{L}^d\text{-a.e.},
		\, c^{(1)}=\tilde{c}^{(1)} \, \nuv\text{-a.e.}, \, c^{(2)}=\tilde{c}^{(2)}\,\nuvinv\text{-a.e.},
			\\\hspace*{1.7 cm} \text{otherwise.}
		\end{cases}
		\end{align*}		
		By abuse of notation, we also identify $\CL^p(\Omega)$ with the space of equivalence classes $\CL^p(\Omega)/\sim$. For $(v,c^{(1)},c^{(2)})\in \CL^p(\Omega)$ we denote the equivalence class by $[(v,c^{(1)},c^{(2)})]$, i.e.,
		\begin{align*}
			[(v,c^{(1)},c^{(2)})]=\{(\tilde{v},\tilde{c}^{(1)},\tilde{c}^{(2)})\in \CL^p(\Omega)\enspace \colon \enspace (v,c^{(1)},c^{(2)})\sim (\tilde{v},\tilde{c}^{(1)},\tilde{c}^{(2)})\}.
		\end{align*}
		Similarly, let $[v],[c^{(1)}],[c^{(2)}]$ be the usual equivalence classes in $\L^1((\Omega,\mathcal{L}^d\lfloor_\Omega);\R)$, $\L^p((\Omega,\nuv);\R^m)$ and $\L^p((\Omega,\nuvinv);\R^m)$, respectively.
		
		\begin{lemma}
			Let $(v,c^{(1)},c^{(2)})\in \CL^p(\Omega)$ with $v=w$ $\mathcal{L}^d$-a.e.\ for some constant $w\in \R$. For $\tilde{w}\in\R$, let $v_{\tilde{w}} \in \L^1((\Omega,\mathcal{L}^d\lfloor_\Omega);\R)$ satisfy $v_{\tilde{w}}=\tilde{w}$ $\mathcal{L}^d$-a.e. If $w\in \R \backslash \{0,1\}$, then
			\begin{align*}
			[(v,c^{(1)},c^{(2)})]=\cup_{\tilde{w}\in \R\backslash \{0,1\}} [v_{\tilde{w}}] \times [c^{(1)}] \times [c^{(2)}].
			\end{align*}
			If $w\in \{0,1\}$, then
			\begin{align*}
			[(v,c^{(1)},c^{(2)})]=[v]\times [c^{(1)}] \times [c^{(2)}].
			\end{align*}
		\end{lemma}
		
		\begin{proof}
			Clearly, if $v=w$ $\mathcal{L}^d$-a.e.\ for $w\in\R \backslash \{0,1\}$, then
			\begin{align*} \frac{|v|}{\|v\|_{\L^1}}=\frac{1}{|\Omega|},\qquad \frac{|1-v|}{\|1-v\|_{\L^1}}=\frac{1}{|\Omega|} \qquad \mathcal{L}^d\text{-a.e.},
			\end{align*}
			independently of the value of $w$.
			For $w\in \{0,1\}$ the claim immediately follows from the definition of the equivalence relation.
			\end{proof}
			
			For the compactness property and $\Gamma$-convergence, we can restrict ourselves to $(v,c^{(1)},c^{(2)})\in \CL^p(\Omega)$ with $0\leq v\leq 1$ $\mathcal{L}^d$-a.e. To see this, note that for any sequence $\{(v_n,c_n^{(1)},c_n^{(2)})\}$ in $\CL^p(\Omega)$ and $\epsilon_n\to 0$ such that $\sup_{n\in\N} \E_{\mu_{\epsilon_n},\epsilon_n} (v_n,c^{(1)}_n,c^{(2)}_n)<+\infty$ we have $v_n\to v$ in $\L^1(\Omega;\R)$ with $v=\chi_E$ for some $E\subset \Omega$. We may consider
			\begin{align}\label{eq:pcapproxvn}
			u_n(x):=\begin{cases}0, & v_n(x)\leq \frac{1}{2},\\
			1, & v_n(x)> \frac{1}{2}, \end{cases}
			\end{align}
			instead of $v_n$. To be precise,
			\begin{lemma}\label{lem:unvnapprox}
				Let $v_n\to v$ in $L^1(\Omega; \R)$, with $v=\chi_E$ for some $E\subset \Omega$. Then
				$\{u_n\}$ defined by \eqref{eq:pcapproxvn} satisfies
				\begin{align*}
				\lim_{n\to \infty} \|u_n-v_n\|_{\L^1(\Omega;\R)}=0.
				\end{align*}
			\end{lemma}
			
			\begin{proof}
				We have
				\begin{align*}
				&\int_\Omega |u_n-v_n|\di x= \int_{\{v_n>\frac{1}{2}\}} |1-v_n|\di x+ \int_{\{v_n\leq\frac{1}{2}\}} |v_n|\di x
				\\&\leq \int_{\{v_n>\frac{1}{2}\}\cap E} |\chi_E-v_n|\di x+ \int_{\{v_n>\frac{1}{2}\}\backslash E} (1+v_n)\di x+ \int_{\{v_n\leq\frac{1}{2}\}\cap E} |v_n|\di x+ \int_{\{v_n\leq\frac{1}{2}\}\backslash E} |v_n-\chi_E|\di x
				\\&\leq \int_{\Omega} |\chi_E-v_n|\di x+ \int_{\Omega\backslash E} 3|v_n|\di x+ \int_{\{v_n\leq\frac{1}{2}\}\cap E} |v_n-\chi_E|\di x+ \int_{\Omega} |v_n-\chi_E|\di x,
				\end{align*}
				where all terms go to 0 as $n\to\infty$ since $v_n\to v$ in $\L^1(\Omega; \R)$.
			\end{proof}
			
			For $(v,c^{(1)},c^{(2)}), (\tilde{v},\tilde{c}^{(1)},\tilde{c}^{(2)}) \in \CL^p(\Omega)$ satisfying $(v,c^{(1)},c^{(2)}) \sim (\tilde{v},\tilde{c}^{(1)},\tilde{c}^{(2)})$,
			where $v$ is nonconstant $\mathcal{L}^d$-a.e.,  $0\leq v\leq 1$ $\mathcal{L}^d$-a.e., and $0\leq \tilde{v} \leq 1$ $\mathcal{L}^d$-a.e., we have $v=\tilde{v}$ $\mathcal{L}^d$-a.e. To see this, note that the equivalence relation on $\CL^p(\Omega)$ implies $v=a\tilde{v}$ $\mathcal{L}^d$-a.e.\ and $1-v=b(1-\tilde{v})$ $\mathcal{L}^d$-a.e.\ for some $a,b \in \R$. For $a\neq b$, we obtain $v=\frac{a(1-b)}{a-b}$ $\mathcal{L}^d$-a.e., in contradiction to $v$ being nonconstant $\mathcal{L}^d$-a.e. This implies that $a=b=1$, i.e., $v=\tilde{v}$ $\mathcal{L}^d$-a.e.
		
		For $(v,c^{(1)},c^{(2)})$ and $(\tilde{v},\tilde{c}^{(1)},\tilde{c}^{(2)})$ in $CL^p(\Omega)$ we define 
		\begin{align*}
			d_{\CL^p}((v,c^{(1)},c^{(2)}),(\tilde{v},\tilde{c}^{(1)},\tilde{c}^{(2)})) := 	d_{\TL^p}((\nuv,c^{(1)}),(\lambda_{|\tilde{v}|},\tilde{c}^{(1)})) +d_{\TL^p}((\nuvinv,c^{(2)}),(\lambda_{|1-\tilde{v}|},\tilde{c}^{(2)}))
		\end{align*}
		where for $(\mu,f),(\lambda,g)$ in $\TL^p(\Omega)$, with
		\begin{align*}
		\TL^p(\Omega):=\{(\mu,f)\colon \mu\in \mathcal{P}(\Omega),f\in \L^p(\Omega,\mu)\},
		\end{align*}
		the metric
		\begin{align*}
		d_{\TL^p}((\mu,f),(\lambda,g)):=	\inf_{\pi\in\Pi(\mu,\lambda)} \bl \int_{\Omega\times \Omega} |x-y|^p + |f(x) - g(x)|^p \di \pi(x,y)\br^\frac{1}{p}
		\end{align*}
		is introduced in \cite{GarciaLimit2016}. 	
		If $\mu,\lambda$ have densities, we can write the distance $d_{\TL^p}$ in the Monge formulation. To be precise,
		\begin{align*}
	        d_{\TL^p}((\mu,f),(\lambda,g))=	\inf_{T\colon {T}_\#\mu= \lambda} \bl \int_{ \Omega} \left[|x-T(x)|^p + |f(x) - g(T(x))|^p \right]\di \mu(x)\br^\frac{1}{p}.
	    \end{align*}
		
		\begin{proposition}
			$(\CL^p(\Omega),d_{\CL^p})$ is a metric space.
		\end{proposition}
	
		\begin{proof}
			Nonnegativity, symmetry  and $d_{\CL^p}((v,c^{(1)},c^{(2)}),(\tilde{v},\tilde{c}^{(1)},\tilde{c}^{(2)}))=0$ for $(v,c^{(1)},c^{(2)})=(\tilde{v},\tilde{c}^{(1)},\tilde{c}^{(2)})$ follow easily from the definition of $d_{\TL^p}$. If $d_{\CL^p}((v,c^{(1)},c^{(2)}),(\tilde{v},\tilde{c}^{(1)},\tilde{c}^{(2)}))=0$, then  $$d_{\TL^p}((\nuv,c^{(1)}),(\lambda_{|\tilde{v}|},\tilde{c}^{(1)}))=0, \qquad d_{\TL^p}((\nuvinv,c^{(2)}),(\lambda_{|1-\tilde{v}|},\tilde{c}^{(2)}))=0,$$
			i.e., $\nuv=\nuvtilde$ $\mathcal{L}^d$-a.e., $\nuvinv=\nuvtildeinv$ $\mathcal{L}^d$-a.e., $c^{(1)}=\tilde{c}^{(1)}$ $\nuv$-a.e., $c^{(2)}=\tilde{c}^{(2)}$ $\nuvinv$-a.e., and these imply
			$$\frac{|v|}{\|v\|_{\L^1}}=\frac{|\tilde{v}|}{\|\tilde{v}\|_{\L^1}}\enspace \mathcal{L}^d\text{-a.e.},\qquad  \frac{|1-v|}{\|1-v\|_{\L^1}}=\frac{|1-\tilde{v}|}{\|1-\tilde{v}\|_{\L^1}} \enspace\mathcal{L}^d\text{-a.e.}$$
			Hence, $(v,c^{(1)},c^{(2)})\sim (\tilde{v},\tilde{c}^{(1)},\tilde{c}^{(2)})$, and we have equality in $\CL^p(\Omega)$.
		\end{proof}

	It was shown in \cite[Proposition 3.12]{GarciaLimit2016} that for $(\mu,f)\in \TL^p(\Omega)$ and a sequence $\{(\mu_n,f_n)\}$ in $\TL^p(\Omega)$,  $(\mu_n,f_n)\to(\mu,f)$ in $\TL^p(\Omega)$ as $n\to \infty$ if and only if $\{\mu_n\}$ converges weakly-$\ast$ to $\mu$ and $f_n\circ T_n\to f$ in $\L^p(\mu)$ as $n\to \infty$ for any stagnating sequence of transportation maps $\{T_n\}$ between $\mu_n$ and $\mu$ with ${T_n}_{\#} \mu=\mu_n$.

	\begin{proposition}\label{prop:convergenceclp}
		Let $(v,c^{(1)},c^{(2)})\in \CL^p(\Omega)$, and let $\{(v_n,c^{(1)}_n,c^{(2)}_n)\}$ be a sequence in $\CL^p(\Omega)$. 
		Then, $(v_n,c_n^{(1)},c_n^{(2)})\to (v,c^{(1)},c^{(2)})$ in $\CL^p(\Omega)$ if and only if $\{\nuvn\}$ converges weakly-$\ast$-$\ast$ to $\nuv$, $c_n^{(1)}\circ T_n^{(1)} \to c^{(1)}$ in $\L^p((\Omega,\nuv);\R^m)$ and $c_n^{(2)}\circ T_n^{(2)}\to c^{(2)}$ in $\L^p((\Omega,\nuvinv);\R^m)$ as $n\to \infty$, for any sequences  of transportation maps $\{T_n^{(1)}\}$ and $\{T_n^{(2)}\}$ satisfying ${T_n^{(1)}}_\#\nuv= \nuvn$, ${T_n^{(2)}}_\#\nuvinv= \nuvinvn$, and $\|T_n^{(1)}-\operatorname{Id}\|_{\L^p(\nuv)}\to 0$, $\|T_n^{(2)}-\operatorname{Id}\|_{\L^p(\nuvinv)}\to 0$.
		\end{proposition}	
	\begin{proof}
 Assume that $(v_n,c_n^{(1)},c_n^{(2)})\to (v,c^{(1)},c^{(2)})$ in $\CL^p(\Omega)$. 
  We have that $$d_{\TL^p}((\nuv,c^{(1)}),(\nuvn,c_n^{(1)}))\to 0$$ and, by \cite[Proposition 3.12]{GarciaLimit2016},  $\{\nuvn\}$ converges weakly-$\ast$ to $\nuv$  and $c_n^{(1)}\circ T_n^{(1)} \to c^{(1)}$ in $\L^p((\Omega,\nuv);\R^m)$ for any sequence of  transportation maps $\{T_n^{(1)}\}$ satisfying the conditions in the proposition. Analogously, we obtain $c_n^{(2)}\circ T_n^{(2)} \to c^{(2)}$ in $\L^p((\Omega,\nuvinv);\R^m)$.

If $\{\nuvn\}$ converges weakly-$\ast$ to $\nuv$, $c_n^{(1)}\circ T_n^{(1)} \to c^{(1)}$ in $\L^p((\Omega,\nuv);\R^m)$ and $c_n^{(2)}\circ T_n^{(2)} \to c^{(2)}$ in $\L^p((\Omega,\nuvinv);\R^m)$, then we conclude that $$d_{\TL^p}((\nuv,c^{(1)}),(\nuvn,c_n^{(1)}))\to 0 \text{ and } d_{\TL^p}((\nuvinv,c^{(2)}),(\nuvinvn,c_n^{(2)}))\to 0.$$ Hence, we obtain that $d_{\CL^p}((v,c^{(1)},c^{(2)}),(\tilde{v},\tilde{c}^{(1)},\tilde{c}^{(2)})) \to 0$.
	\end{proof}

\section{\texorpdfstring{$\Gamma$}{Gamma}-convergence for  piecewise constant segmentations}\label{sec:pw_constant}

In this section we study the Ginzburg-Landau image segmentation model where $c^{(1)},c^{(2)}\in\R^m$ are constants and correspond to the optimal intensity values to approximate each of the two segments. 
For constants $c^{(1)},c^{(2)}$,  we define  $\bar{\E}_\epsilon\colon \L^1(\Omega;\R)\times \R^m\times \R^m$  by
\begin{align}\label{eq:energyepspc}
\bar{\E}_{\epsilon}(v,c^{(1)},c^{(2)}):=\begin{cases}\int_{\Omega}\left( |c^{(1)}-u_0|^p |v| + |c^{(2)}-u_0|^p |1-v|\right) \di x +\frac{\nu}{c_W} \int_{\Omega}\left(\epsilon|\nabla v|^2+\frac{1}{\epsilon}W(v)\right) \di x, \\  \hspace*{1.7cm} \text{if }v\in \W^{1,2}(\Omega;\R),c^{(1)},c^{(2)}\in\R^m,\\ +\infty, \\  \hspace*{.85 cm} \text{otherwise,}
\end{cases}
\end{align}
where $c_W$ is defined in \eqref{eq:cw}, and $u_0\in L^\infty(\Omega;\R^m)$ is  given. 
The aim of this section is to show that $\{\bar{\E}_{\epsilon}\}$ $\Gamma$-converges to $\bar{\E}
\colon \L^1(\Omega;\R)\times \R^m\times \R^m$, defined by
\begin{align}\label{eq:energypc}
\bar{\E}(v,c^{(1)},c^{(2)}):=\begin{cases}
\int_{E} |c^{(1)}-u_0|^p\di x+\int_{\Omega\backslash E} |c^{(2)}-u_0|^p \di x +\nu  \tv(v), \\\hspace*{1.7cm}\text{if }v=\chi_E\in\bv(\Omega;\{0,1\}) \text{ for }E:=\{x\in\Omega\colon v(x)=1\},
 \\  \hspace*{1.7cm}  c^{(1)},c^{(2)}\in \R^m, 
\\ +\infty, \hspace*{0.85 cm}  \text{otherwise.}	\end{cases}
\end{align}
Note that $\bar{\E}_\epsilon$ and $\bar{\E}$ follow immediately from the definition of $\E_{\mu,\epsilon}$ and $\E_\mu$ when $c^{(1)},c^{(2)}$  constant. In this case, the $\CL^p(\Omega)$ topology  is not practical and we consider the $\L^1(\Omega;\R)\times \R^m\times \R^m$ topology instead.  The main results of this section are the compactness property and the $\Gamma$-convergence of $\bar{\E}_\epsilon$ for piecewise constant segmentations, which imply the convergence of minimizers:
\begin{theorem}\label{th:gammaconvpc}
	Let $\Omega \subset \R^d$ be an open, bounded  set, let $1<p<+\infty$, and let $\bar{\E}_\epsilon\colon \L^1(\Omega;\R)\times \R^m\times \R^m\to [0,+\infty]$ and $\bar{\E}
	\colon \L^1(\Omega;\R)\times \R^m\times \R^m$ as defined by \eqref{eq:energyepspc} and \eqref{eq:energypc}, respectively.
	Then, the functional $\bar{\E}_{\epsilon}$ satisfies the compactness property and $\Gamma$-converges with respect to the $\L^1(\Omega;\R)\times\R^m\times \R^m $ topology to $\bar{\E}$ as $\epsilon\to 0$.
\end{theorem}

Let us first state a general lemma which is not only valid for constant functions $c^{(1)},c^{(2)}$, but more generally for  functions  $c^{(1)}\in W^{1,p}((\Omega,\nuv);\R^m), c^{(2)}\in \W^{1,p}((\Omega,\nuvinv);\R^m)$ for $v\in \L^1((\Omega,\mathcal{L}^d\lfloor_\Omega);\R)$ given.

\begin{lemma}\label{lem:compactnesshelp}
	Let $\Omega\subset\R^d$ be an open set with finite measure.  
	Define the energy functional $\E_{\mu_\epsilon,\epsilon}$ as in \eqref{eq:approxenergy}, and let $\epsilon_n\to 0$, $\{v_n\}\subset  \W^{1,2}((\Omega,\mathcal{L}^d\lfloor_\Omega);\R)$,   $\{c^{(1)}_n\},\{c^{(2)}_n\}$ such that 
	$c^{(1)}_n\in W^{1,p}((\Omega,\nuvn);\R^m)$, $c^{(2)}_n\in \W^{1,p}((\Omega,\nuvinvn);\R^m)$, and assume that
	\begin{align*}
	M:=\sup_{n\in\N}\E_{\mu_{\epsilon_n},\epsilon_n}(v_n,c^{(1)}_n,c^{(2)}_n)<+\infty,
	\end{align*}
	where $\mu_{\epsilon_n} \to \mu\in(0,+\infty]$ as $\epsilon_n\to 0$.
	Then, there exist  a subsequence $\{v_{n_k}\}$ of  $\{v_n\}$ and $v\in \bv(\Omega;\{0,1\})$, with $v=\chi_E$ $\mathcal{L}^d$-a.e.\ for a Lebesgue measurable set $E\subset \Omega$, such that 
	${v_{n_k}\to v}$ in $\L^1(\Omega;\R)$. 
\end{lemma}

\begin{proof}
	Since $\sup_{n\in\N}\E_{\epsilon_n}^{GL}(v_n)<+\infty,$ where $\E_\epsilon^{GL}$ denotes the Ginzburg-Landau energy functional defined in \eqref{eq:energyginzburglandau}, Theorem \ref{th:compactness} can be invoked.
\end{proof}

As a first step towards proving Theorem \ref{th:gammaconvpc}, we show a compactness result based on Lemma \ref{lem:compactnesshelp}.
\begin{theorem}[Compactness]\label{th:compactnessconst}
	Let $\Omega\subset\R^d$ be an open set with finite measure, let $\epsilon_n\to 0$, and let $\{v_n\}\subset \W^{1,2}(\Omega;\R)$,  $\{c^{(1)}_n\},\{c^{(2)}_n\}\subset\R^m$, be such that
	\begin{align*}
	M:=\sup_{n\in\N}\bar{\E}_{\epsilon_n}(v_n,c^{(1)}_n,c^{(2)}_n)<+\infty.
	\end{align*}
	Then, there exist  a subsequence $\{v_{n_k}\}$ of  $\{v_n\}$ and $v\in \bv(\Omega;\{0,1\})$, with $v=\chi_E$ for some Lebesgue measurable set $E\subset \Omega$, such that 
	${v_{n_k}\to v}$ in $L^1(\Omega;\R)$. If $\mathcal{L}^d(E)>0$, then there exists a converging subsequence $\{c^{(1)}_{n_k}\}$ of $\{c^{(1)}_n\}$ with limit ${c^{(1)}\in\R^m}$. If $\mathcal{L}^d(\Omega\backslash E)>0$, then there exists a converging subsequence $\{c^{(2)}_{n_k}\}$ of $\{c^{(2)}_n\}$ with limit ${c^{(2)}\in\R^m}$.
\end{theorem}
\begin{proof}
	By Lemma \ref{lem:compactnesshelp} we can find  a  subsequence $\{v_{n_k}\}$ of  $\{v_n\}$ and $v\in \bv(\Omega;\{0,1\})$, with $v=\chi_E$ for some Lebesgue measurable set $E\subset \Omega$, such that 
	${v_{n_k}\to v}$ in $\L^1(\Omega;\R)$.
	 For $\mathcal{L}^d(E)>0$ the sequence $\{c^{(1)}_n\}$ has to be bounded. To see this, note that the energy bound  implies that $\int_{\Omega}|c^{(1)}_n-u_0|^p |v_n|\di x$ is uniformly bounded. If the sequence $\{c^{(1)}_n\}$ was unbounded, for every $n>0$ there exists some $k_n\in \N$ such that 
	 $$M\geq \sup_{n\in\N}\int_{\Omega}|c^{(1)}_n-u_0|^p |v_n|\di x \geq  n \| v_{k_n}\|_{L^1(\Omega)}$$
	 using the fact that $u_0$ is bounded. This implies that $v_{k_n}\to 0$ in $L^1(\Omega; \R)$, which contradicts $\mathcal{L}^d(E)>0$.
	 Hence, $\{c^{(1)}_n\}$ is bounded, and the existence of  a subsequence of $\{c^{(1)}_n\}$ converging to $c^{(1)}$ in $\R^m$ follows immediately from the Bolzano–Weierstrass theorem. Similarly, one can show if $\mathcal{L}^d(\Omega\backslash E)>0$ then $\{c^{(2)}_n\}$ is bounded, and has a converging subsequence with limit $c^{(2)}\in\R^m$.	 
\end{proof}

\begin{proof}[Proof of Theorem \ref{th:gammaconvpc}]
Since the compactness property follows from Theorem \ref{th:compactnessconst}, it remains to show the $\Gamma$-convergence. Let
\begin{align}
\bar{\E}^{(1)}(v,c^{(1)},c^{(2)}) :=\int_{\Omega}\left( |c^{(1)}-u_0|^p |v| + |c^{(2)}-u_0|^p |1-v|\right) \di x,
\end{align}
so that $\bar{\E}_\epsilon(v,c^{(1)},c^{(2)}) =\bar{\E}^{(1)}(v,c^{(1)},c^{(2)}) +\frac{\nu}{c_W}\E_{\epsilon}^{GL}(v)$ and (when $v=\chi_E$) $\bar{\E}(v,c^{(1)},c^{(2)}) =\bar{\E}^{(1)}(v,c^{(1)},c^{(2)}) +\nu\tv(v)$.

Let $(v_n,c_n^{(1)},c_n^{(2)})  \to (v,c^{(1)},c^{(2)})$ in $\CL^p(\Omega)$, i.e., $v_n$ is bounded in $\L^1(\Omega;\R)$,  $\sup_{n\in \N} |c_n^{(1)}|<+\infty$ and $\sup_{n\in \N} |c_n^{(2)}|<+\infty$. 
We have 
\begin{align*}
&|\bar{\E}^{(1)}(v_n,c_n^{(1)},c_n^{(2)})-\bar{\E}^{(1)}(v,c^{(1)},c^{(2)})| \\&\leq \int_\Omega \left( |c_n^{(1)}-u_0|^p \left||v_n|-| v|\right| + |v| \left| |c_n^{(1)}-u_0|^p-|c^{(1)}-u_0|^p\right| \right)\di x\\
&\quad + \int_\Omega \left( |c_n^{(2)}-u_0|^p \left||1-v_n|-| 1-v|\right| + |1-v| \left| |c_n^{(2)}-u_0|^p-|c^{(2)}-u_0|^p\right| \right)\di x.
\end{align*}
Note that for any $\delta>0$ there exists $C_\delta>0$ such that for all $a,b \in \R^m$ we have $$|a|^p\leq (1+\delta) |b|^p +C_\delta |a-b|^p,$$
implying
$$|c_n^{(i)}-u_0|^p\leq (1+\delta) |c^{(i)}-u_0|^p+C_\delta |c_n^{(i)}-c^{(i)}|^p.$$
Hence, 
\begin{align*}
&|\bar{\E}^{(1)}(v_n,c_n^{(1)},c_n^{(2)})-\bar{\E}^{(1)}(v,c^{(1)},c^{(2)})| \\
&\leq \bl \sup_{x\in\Omega} |c_n^{(1)}-u_0(x)|^p\br \|v_n-v\|_{\L^1(\Omega)}+\int_\Omega |v| \bl \delta |c^{(1)}-u_0|^p+C_\delta |c_n^{(1)}-c^{(1)}|^p\br\di x\\
&\quad+ \bl \sup_{x\in\Omega} |c_n^{(2)}-u_0(x)|^p\br \|v_n-v\|_{\L^1(\Omega)}+\int_\Omega |1-v| \bl \delta |c^{(2)}-u_0|^p+C_\delta |c_n^{(2)}-c^{(2)}|^p\br\di x\\
&\leq C \|v_n-v\|_{\L^1(\Omega)} + \delta \bar{\E}^{(1)}(v,c^{(1)},c^{(2)})+C_\delta |c_n^{(1)}-c^{(1)}|^p \|v\|_{\L^1(\Omega)}+C_\delta |c_n^{(2)}-c^{(2)}|^p \|1-v\|_{\L^1(\Omega)}.
\end{align*}
Letting $n\to \infty$, we have
\begin{align*}
\lim_{n\to \infty}|\bar{\E}^{(1)}(v_n,c_n^{(1)},c_n^{(2)})-\bar{\E}^{(1)}(v,c^{(1)},c^{(2)})| \leq \delta \bar{\E}^{(1)}(v,c^{(1)},c^{(2)})
\end{align*}
for any $\delta>0$. Let $\delta\to 0$ to obtain
\begin{align*}
\lim_{n\to \infty} \bar{\E}^{(1)}(v_n,c_n^{(1)},c_n^{(2)})=\bar{\E}^{(1)}(v,c^{(1)},c^{(2)}).
\end{align*}
By the stability of $\Gamma$-convergence under its continuous perturbations \cite[Proposition 6.20]{dalmaso1992}, we obtain the $\Gamma$-convergence of $\bar{\E}_\epsilon$ to $\{\bar{\E}\}$ in $L^1(\Omega;\R)\times \R^m\times \R^m$.
\end{proof}

Due to the compactness result in Theorem \ref{th:compactnessconst}, we only consider $\emptyset \subsetneq E \subsetneq \Omega$ with $0<\mathcal{L}^d(E)<\mathcal{L}^d(\Omega)$ for minimizers of the function $\bar{\E}$ in \eqref{eq:energypc}. However, the $\Gamma$-limit $\bar{\E}$ in \eqref{eq:energypc} is defined for all sets $\emptyset \subset E \subset \Omega$.

\section{\texorpdfstring{$\Gamma$}{Gamma}-convergence for piecewise smooth approximations}\label{sec:pw_smooth}
In this section we prove the main result of the paper, stated in Theorem \ref{th:mainresult}, namely the $\Gamma$-convergence of the energy functional $\E_{\mu_\epsilon, \epsilon}$ in \eqref{eq:approxenergy} for any positive parameter $\mu_\epsilon$. 
In the following we differentiate between two regimes depending on the convergence of the positive parameter $\mu_{\epsilon_n}$ as $\epsilon_n\to 0$:
\begin{enumerate}
	\item $\mu_{\epsilon_n}\to \mu$ for a constant $\mu>0$,
	\item $\mu_{\epsilon_n}\to +\infty$.
\end{enumerate}
These two cases  cover all positive limits of $\mu_{\epsilon_n}$ as $\epsilon_n\to 0$. We note that the analysis is very similar for  $\lim_{\epsilon_n\to 0}\mu_{\epsilon_n}=\mu>0$ and $\lim_{\epsilon_n\to 0}\mu_{\epsilon_n} =+\infty$.
We start by showing compactness:
\begin{theorem}[Compactness]\label{th:compactnesssmoothcase1}
	Let $\Omega\subset\R^d$, with $d\geq 2$, be an open set with finite measure. Let $\epsilon_n\to 0$ and $\{v_n\}\subset \W^{1,2}((\Omega,\mathcal{L}^d\lfloor_\Omega);\R),\{c^{(1)}_n\},\{c^{(2)}_n\}$ 
	be such that  $c_n^{(1)}\in \W^{1,p}((\Omega,\nuvn);\R^m)$, $c_n^{(2)}\in \W^{1,p}((\Omega,\nuvinvn);\R^m)$, 
	and
	\begin{align*}
	M:=\sup_{n\in\N}\E_{\mu_{\epsilon_n},\epsilon_n}(v_n,c^{(1)}_n,c^{(2)}_n)<+\infty,
	\end{align*}
	for $\E_{\mu_\epsilon, \epsilon}$ defined in \eqref{eq:approxenergy}, with $\lim_{n\to\infty}\mu_{\epsilon_n }\in(0,+\infty]$.
	Then, there exist a subsequence $\{v_{n_k}\}$ of $\{v_n\}$ and $v\in \bv(\Omega;\{0,1\})$, with $v=\chi_E$ for some Lebesgue measurable set $E\subset \Omega$, such that 
	${v_{n_k}\to v}$ in $\L^1(\Omega;\R)$.
	If $\mathcal{L}^d(E)>0$,   
if there exist $\kappa>0$, $r_0>0$  such that $P(E;B_r(x)) \geq \kappa r^d$ for every $x \in \partial^* E$, and $\{c^{(1)}_n\}_{n\in\N}$ are bounded in $\L^\infty$, then $(\nuvn,c_n^{(1)})$ is precompact in $\TL^p$, and any cluster point $(\nuv,c^{(1)})$ satisfies $c^{(1)}\in \W^{1,p}((\Omega,\nuv);\R^m)$.  
    Similarly, if $\mathcal{L}^d(\Omega\setminus E)>0$, 
  if  there exist ${\kappa>0}$, $r_0>0$  such that $P(\Omega \setminus E;B_r(x)) \geq \kappa r^d$ for every $x \in \partial^* (\Omega \setminus E)$,  and $\{c^{(2)}_n\}_{n\in\N}$ are bounded in $\L^\infty$, then $(\nuvinvn,c_n^{(2)})$ is precompact in $\TL^p$, and any cluster point $(\nuvinv,c^{(2)})$ satisfies $c^{(2)}\in \W^{1,p}((\Omega,\nuvinv);\R^m)$. 
	In particular, if $0<\mathcal{L}^d(E)<\mathcal{L}^d(\Omega)$, if the above assumptions on the
 perimeter of $E$ and of $\Omega \setminus E$ hold, and if $\{c^{(i)}_n\}_{n\in\N}$, $i=1,2$, are bounded in $\L^\infty$, then 
	there exist a subsequence $(v_{n_k},c^{(1)}_{n_k},c^{(2)}_{n_k})$ of $(v_n,c^{(1)}_n,c^{(2)}_n)$ and $(v,c^{(1)},c^{(2)})\in \CL^p(\Omega)$ such that $\{(v_{n_k},c^{(1)}_{n_k},c^{(2)}_{n_k})\}$ converges to $(v,c^{(1)},c^{(2)})$ in $\CL^p(\Omega)$ and $\E_{\mu}(v,c^{(1)},c^{(2)})<+\infty$. 
\end{theorem}

\begin{proof}
The existence of a subsequence $\{v_{n_k}\}$ of $\{v_n\}$ and $v\in \bv(\Omega;\{0,1\})$ with $v=\chi_E$ for a measurable set $E\subset \Omega$ with finite perimeter such that $v_{n_k}\to v$ in $\L^1(\Omega;\R)$ follows from Lemma~\ref{lem:compactnesshelp}. 
In particular, $\{\nuvn\}$ and $\{\nuvinvn\}$ converge weakly-$\ast$ to $\nuv$ and $\nuvinv$, respectively. 

Let us first consider $0<\lim_{n\to \infty}\mu_{\epsilon_n}<+\infty$, and we  assume, without loss of generality, that $\mu_{\epsilon_n}$ are uniformly bounded by  positive constants  from above and below.
Since the existence of converging subsequences of $\{c_{n_{k_l}}^{(1)}\circ T_{n_{k_l}}^{(1)}\}\to c^{(1)}$ and $\{c_{n_{k_l}}^{(2)}\circ T_{n_{k_l}}^{(2)}\}\to c^{(2)}$ can be shown in a similar way, we restrict ourselves to $c_{n_{k_l}}^{(1)}\circ T_{n_{k_l}}^{(1)}\to c^{(1)}$ and in the following assume that $\mathcal{L}^d(E)>0$.
For ease of notation, we omit the superscript index $(1)$.  

Since $\{\nuvn\}$ converges weakly-$\ast$ to $\nuv$, then $\{\nuvn\}$ converges in the $p^\prime$-Wasserstein distance to $\nuv$, with  $\frac{1}{p}+\frac{1}{p^\prime}=1$.
In particular, there exists a sequence of transport maps $\{T_n\}$ satisfying
\[ T_{n\#}\nuv= \nuvnk \quad \text{and}\quad \lim_{n\to\infty} \| T_n-\operatorname{Id}\|_{\L^{p^\prime}(E)}=0. \]

Let $\psi\in \C_c^\infty(\R^d)$ be a standard mollifier, e.g., 
\[ \psi (x):=\begin{cases} C \exp\bl \frac{1}{|x|^2-1}\br,& |x|<1, \\ 0, &  |x|\geq 1,\end{cases} \]
where the constant $C>0$ is chosen such that $\int_{\R^d} \psi \di x=1$.
For each $a>0$, we set 
\[ \psi_{a}(x)=\frac{1}{a^d}\psi\bl \frac{x}{a}\br,\quad x\in\R^d. \]
We define convolution in the usual way, i.e., $(\psi\ast c)(x):= \int_\Omega \psi(x-y) c(y)\di y$, and for convenience we let $\hat{v}(x):=\frac{v(x)}{|E|}$.
We claim that there exists a positive converging sequence $\{a_n\}_{n\in\N}\subset \R$ with $\lim_{n\to \infty} a_n = 0$, such that
\begin{equation} \label{eq:liminfcompactass1}
\sup_{n\in\N} \|\nabla (\psi_{a_n}\ast ((c_n\circ T_n) \hat{v}))\|_{\L^{1}(\R^d)} <+\infty,
\end{equation}
and 
\begin{equation} \label{eq:liminfcompactass2}
\lim_{n\to\infty} \| \psi_{a_n} \ast ((c_n\circ T_n) \hat{v}) -(c_n\circ T_n) \hat{v} \|_{\L^{1}(E)}=0.
\end{equation}
Under these assumptions, we show that $c\in \W^{1,p}((\Omega,\nuv;\R^m)$.
Since
\begin{align*}
\| \psi_{a_n}\ast ((c_n\circ T_n) \hat{v})\|_{\L^p(\R^d)}^p & = \int_{\R^d} \left| \frac{1}{|E|} \int_E \psi_{a_n}(x-y) c_n(T_n(y)) \, \mathrm{d} y \right|^p \, \mathrm{d} x \\
 & \leq \frac{1}{|E|^p} \int_{\R^d} \int_E \psi_{a_n}(x-y) |c_n(T_n(y)|^p \, \mathrm{d} y \, \mathrm{d} x \\
 &=  \frac{1}{|E|^p} \int_E  |c_n(T_n(y)|^p \, \mathrm{d} y= \frac{1}{|E|^{p-1}} \| c_n\circ T_n\|_{\L^p(\nuv)}^p,
\end{align*}
by the boundedness of  $\{\E_{\mu_{\epsilon_n},\epsilon_n}\}$, we deduce that the sequence $\| \psi_{a_n}\ast ((c_n\circ T_n) \hat{v})\|_{\L^p(\R^d)}$ is uniformly bounded in $n$.
The sequence $\{f_n\}$, with $f_n := \psi_{a_n}\ast ((c_n\circ T_n) \hat{v})$, is bounded in $\mathrm{W}^{1,1}(\R^d)$ by \eqref{eq:liminfcompactass1}, and so there exists a subsequence (not relabelled) that converges in $\L^\alpha(\R^d)$ to some $f\in \mathrm{W}^{1,1}(\R^d)$ for any $1\leq \alpha\leq\frac{d}{d-1}$. 
Choosing $\alpha=1$, we have
\[ \lim_{n\to\infty} \| \psi_{a_n} \ast ((c_n\circ T_n) \hat{v}) - f \|_{\L^1(\R^d)}=0. \]
Define $c(x) := f(x) |E|$ for $x\in E$ and $c(x)=0$ for $x\not\in E$. We obtain
\[ \| \psi_{a_n} \ast ((c_n\circ T_n) \hat{v}) - c \hat{v} \|_{\L^1(E)} \leq \| \psi_{a_n} \ast ((c_n\circ T_n) \hat{v}) - f \|_{\L^1(\R^d)} \to 0. \]
Together with~\eqref{eq:liminfcompactass2}, and the fact that
\begin{align*}
\| c_n\circ T_n -c \|_{\L^1(\nuv)} & = \| (c_n\circ T_n) \hat{v} - c\hat{v} \|_{\L^1(E)} \\
 & \leq \| (c_n\circ T_n) \hat{v} - \psi_{a_n} \ast ((c_n\circ T_n) \hat{v})  \|_{\L^1(E)} + \| \psi_{a_n} \ast ((c_n\circ T_n) \hat{v}) -c \hat{v} \|_{\L^1(E)},
\end{align*} 
we deduce that
\[ \lim_{n\to\infty} \| c_n\circ T_n -c \|_{\L^1(\nuv)}=0. \]
From the fact that $\{c_n\}$ is bounded in $\L^\infty$,  and extracting a further subsequence (not relabeled), we can assume that $c_n\circ T_n\to c$ pointwise on $E$, which implies that  $c$ is also bounded in $\L^\infty$.
Moreover,
\[ \| c_n\circ T_n -c \|_{\L^p(\nuv)}^p \leq \| c_n\circ T_n -c \|_{\L^\infty(E)}^{p-1} \| c_n\circ T_n -c \|_{\L^1(\nuv)} \to 0. \]

Since $c_n\in W^{1,p}(\nuvn)$, by Theorem \ref{th:relationdiffmetricmeasure} there exist $g_n\in\L^p(\nuvn)$ and $\Omega_n\subset \Omega$ 
satisfying
\[ |c_n(x) - c_n(y)|\leq |x-y| (g_n(x) + g_n(y)) \]
for all $x,y\in \Omega_n$, $\nuvn(\Omega_n)=1$, 
and $\|g_n\|_{\L^p(\nuvn)}\leq \|c_n\|_{\L^{1,p}(\nuvn)}+\frac{1}{n}$.
Using the transport maps $T_n$, we can rewrite the above as
\begin{equation} \label{eq:cnDer}
|c_n(T_n(x)) - c_n(T_n(y))|\leq |T_n(x)-T_n(y)| (g_n(T_n(x)) + g_n(T_n(y)))
\end{equation}
for all $x,y\in T_n^{-1}(\Omega_n)\subset E$,
$\nuv(T_n^{-1}(\Omega_n))=1$, 
and $\|g_n\circ T_n\|_{\L^p(\nuv)}\leq \|c_n\|_{\L^{1,p}(\nuvn)}+1$.
Now $\nuv(T_n^{-1}(\Omega_n))=1$ implies $|T_n^{-1}(\Omega_n)|=|E|$, and so~\eqref{eq:cnDer} holds for almost every $x,y\in E$.
Taking the  union over all sets $T_n^{-1}(\Omega_n)$ for $n\in\N$, we can further say that there exists $\tilde{E}$ with $|\tilde{E}|=|E|$ such that~\eqref{eq:cnDer} holds for all $x,y\in\tilde{E}$ and $n\in\N$.
As $g_n\circ T_n$ are bounded in $\L^p(\nuv)$,  there exists a  weakly converging subsequence to some $g\in \L^p(\nuv)$. 
Moreover, $(x,y)\mapsto \chi_{T_n^{-1}(\Omega_n)^2}(x,y)|T_n(x)-T_n(y)| (g_n(T_n(x))+g_n(T_n(y)))$ is also bounded in $\L^{p}(\lambda_{|\nu|})$, and so it converges weakly along a subsequence to $(x,y)\mapsto \chi_{\Omega^2}(x,y)|x-y|(g(x)+g(y))$.
For any $\phi:\R^d\times\R^d\to \R_+$, with $\phi\in\L^{p^*}(\lambda_{|\nu|}\times \lambda_{|\nu|})$ where $p^*>0$ such that $\tfrac{1}{p}+\tfrac{1}{p^*}=1$, we have
\begin{align*}
& \int_{\Omega^2} |c(x) - c(y)| \phi(x,y) \di x \di y \\
& \qquad \qquad \leq \liminf_{n\to\infty} \int_{T_n^{-1}(\Omega_n)^2} |c_n(x) - c_n(y)| \phi(x,y) \di x \di y \qquad \text{by Fatou's lemma} \\
& \qquad \qquad \leq \liminf_{n\to\infty} \int_{T_n^{-1}(\Omega_n)^2} |T_n(x) - T_n(y)|(g_n(T_n(x)) + g_n(T_n(y))) \phi(x,y) \di x \di y \qquad \text{by \eqref{eq:cnDer}} \\
& \qquad \qquad = \int_{\Omega^2} |x-y| (g(x)+g(y)) \phi(x,y) \di x \di y.
\end{align*}
Therefore,
\[ |c(x) - c(y)|\leq |x-y| (g(x) + g(y)) \]
for almost every $x,y\in \tilde{E}$.
By redefining $g(x)=+\infty$, $g(y)=+\infty$ for any $(x,y)$ where the above does not hold, we can assume that the inequality holds for all $(x,y)$ (and the $\L^p$ norm of $g$ is unchanged).
By the weak lower semi-continuity of norms, we have $\|g\|_{\L^p(\nuv)}\leq \liminf_{n\to\infty} \|g_n\circ T_n\|_{\L^p(\nuv)}$, where the right hand side is finite due to boundedness of the energies.
It follows that $c\in \W^{1,p}((\Omega,\nuv;\R^m)$, and
\begin{align}
\| c\|_{\L^{1,p}(\nuv)} & \leq \|g\|_{\L^p(\nuv)} \notag 
 \leq \liminf_{n\to\infty} \|g_n\circ T_n\|_{\L^p(\nuv)} \notag 
  =  \liminf_{n\to\infty} \|g_n\|_{\L^p(\nuvn)} \notag \\
 & \leq \liminf_{n\to\infty} \bl \|c_n\|_{\L^{1,p}(\nuvn)} + \frac{1}{n} \br \notag \\
 & =  \liminf_{n\to\infty} \|c_n\|_{\L^{1,p}(\nuvn)}. \label{eq:cnl1pliminf}
\end{align}

If $\mu_{\epsilon_n}\to +\infty$ as $n\to \infty$, the existence of a converging subsequence $\{(v_{n_k},c^{(1)}_{n_k},c^{(2)}_{n_k})\}$  with limit in $\CL^p(\Omega)$ follows in the same way as for the case $\lim_{n\to \infty} \mu_{\epsilon_{n}}< +\infty$,  due to the uniform boundedness of $\E_{\mu_{\epsilon_n},\epsilon_n}(v_n,c^{(1)}_n,c^{(2)}_n)$. 
Furthermore, if $\mu_{\epsilon_n}\to +\infty$ as $n\to \infty$, we have that $\{c_n^{(1)}\}$ converges to a constant since, again omitting the superscript $(1)$,
\[ \limsup_{n\to \infty}  \|c_n\|_{\L^{1,p}(\nuvn)} = \limsup_{n\to \infty}  \inf_{g_n}\|g_n\|_{\L^{p}(\nuvn)} = \limsup_{n\to \infty} \inf_{g_n}\|g_n\circ T_n\|_{\L^{p}(\nuv)}=0 \]
i.e., $g_n\circ T_n\to 0$ in $\L^p(\Omega,\nuv)$, and taking the limit on both sides of the following inequality 
\[ |c_n(T_n(x)) - c_n(T_n(y))| \leq |T_n(x)-T_n(y)|(g_n(T_n(x)) + g_n(T_n(y)))  \enspace \nuv\text{-a.e.}, \]
implies that there exists a constant $c_1\in\R^m$ such that $c^{(1)}=c_1$ $\mathcal{L}^d$-a.e.\ $x\in  E$.
Similarly, it follows that  $c^{(2)}=c_2$ $\mathcal{L}^d$-a.e.\ $x\in \Omega \backslash E$ for some constant $c_2\in \R^m$.
	  
It remains to show~\eqref{eq:liminfcompactass1} and~\eqref{eq:liminfcompactass2}.
To show that \eqref{eq:liminfcompactass1} is indeed satisfied, note that we have for any positive converging sequence $\{a_n\}_{n\in\N}\subset \R$ with $\lim_{n\to \infty} a_n= 0$ (which will be specified later),
\begin{align*}
\nabla \left( \psi_{a_n}\ast ((c_n\circ T_n)\hat{v})\right)(x) & = \frac{1}{a_n^d |E|}\nabla \int_E \psi\bl \frac{x-y}{a_n}\br  c_n(T_n(y)) \di y \\ 
 & = \frac{1}{a_n^{d+1} |E|}\int_E \nabla \psi\bl \frac{x-y}{a_n}\br  c_n(T_n(y)) \di y \\
 & = \frac{1}{a_n^{d+1} |E|}\int_E \nabla \psi\bl \frac{x-y}{a_n}\br \bl c_n(T_n(y))-c_n(T_n(x))\br \di y \\
 & \qquad - \frac{1}{a_n^{d+1} |E|}\int_{\R^d\setminus E} \nabla \psi\bl \frac{x-y}{a_n}\br c_n(T_n(x)) \di y, 
\end{align*}
where we  extended $c_n$ to be zero outside of $\Omega$ in the last equality and used the fact that $\int_{\R^d} \nabla \psi (\frac{x-y}{a_n}) \di y=0$ as $\nabla \psi$ is odd. 
Hence,
\begin{align*}
& \|\nabla (\psi_{a_n}\ast ((c_n\circ T_n) \hat{v})) \|_{\L^1(\R^d)} \\
& \leq \frac{1}{a_n^{d+1}|E|} \int_{\R^d} \left| \int_E \nabla \psi\bl \frac{x-y}{a_n}\br \bl c_n(T_n(y))-c_n(T_n(x))\br \di y \right| \di x \\
& \qquad + \frac{1}{a_n^{d+1} |E|} \int_{\R^d} \left| \int_{\R^d\setminus E} \nabla \psi\bl \frac{x-y}{a_n}\br c_n(T_n(x)) \di y \right|\di x \\
& =: \I_n + \II_n. 
\end{align*}

Starting with term $\II_n$, a change of variables implies that
\begin{align*}
\II_n & \leq \frac{\|c_n\|_{\L^\infty}}{a_n^{d+1}|E|} \int_{E} \int_{\R^d\setminus E} |\nabla \psi|\bl \frac{x-y}{a_n} \br \di y \di x \\
 & = \frac{\|c_n\|_{\L^\infty}}{a_n |E|} \int_{x\in E\,:\,\dist(x,\partial E)\leq a_n} \int_{w\,:\,x-a_nw\in \R^d\setminus E} |\nabla \psi|(w) \di w \di x \\
 & \leq \frac{\|c_n\|_{\L^\infty} \|\nabla \psi\|_{\L^\infty} |B(0,1)|}{|E|} \frac{|\{ x\in E\,\colon\, \dist(x,\partial E) \leq a_n\}|}{a_n}.
\end{align*}
Note that the assumption in Corollary \ref{cor:PerBound} that the topological boundary $\partial E$ is the closure of the reduced boundary $\partial^* E$ holds for free up to a modification on a Lebesgue null set, see \cite[Proposition 12.20]{maggi2012sets}. 
By Corollary~\ref{cor:PerBound}, we can choose the sequence $\{a_n\}_{n\in\N}$ such that we have that $\II_n=O(1)$.

For the term $\I_n$, we use~\eqref{eq:cnDer} to infer 
\begin{align*}
I_n & \leq \frac{1}{a_n^{d+1}|E|} \int_E  \left| \int_E |\nabla \psi| \left(\frac{x-y}{a_n}\right) \left| T_n(x) - T_n(y) \right| (g_n(T_n(x)) + g_n(T_n(y))) \di y \right| \di x \\
 & \qquad + \frac{1}{a_n^{d+1} |E|} \int_{\R^d\setminus E} \left| \int_E (\nabla \psi)\bl \frac{x-y}{a_n}\br \bl c_n(T_n(y))\br \di y \right| \di x \\
 & \leq \frac{1}{a_n^{d+1} |E|} \int_E \int_E |\nabla \psi| \left(\frac{x-y}{a_n}\right) \left( 2|T_n(x) - x| + |x-y| \right) (g_n(T_n(x)) + g_n(T_n(y))) \di y \di x \\
 & \qquad + \frac{\|c_n\|_{\L^\infty}}{a_n^{d+1} |E|} \int_{\R^d\setminus E} \int_E |\nabla \psi|\bl\frac{x-y}{a_n}\br \di y \, \di x.
\end{align*}
The second term above can be shown to be $O(1)$ following the same argument as for $\II_n$.

We let
\[ \III_n:= \frac{2}{a_n^{d+1} |E|} \int_E \int_E |\nabla \psi| \left(\frac{x-y}{a_n}\right) |T_n(x) - x| (g_n(T_n(x)) + g_n(T_n(y))) \di y \di x \]
and
\[ \IV_n := \frac{1}{a_n^{d+1} |E|} \int_E \int_E |\nabla \psi| \left(\frac{x-y}{a_n}\right) |x-y| (g_n(T_n(x)) + g_n(T_n(y))) \di y \di x. \]
A change of variables implies
\begin{align*}
\III_n & \leq \frac{2}{a_n |E|} \int_E \int_{z\,:\,x-a_nz\in E} |\nabla \psi|(z) |T_n(x) - x| (g_n(T_n(x)) + g_n(T_n(x-a_nz))) \di z \di x \\
 & \leq \frac{2}{a_n |E|} \int_E |T_n(x) - x| \bl \|\nabla \psi\|_{\L^1} g_n(T_n(x)) + \|\nabla \psi\|_{\L^\infty} \int_{B(0,1)} g_n(T_n(x-a_nz)) \di z \br \di x \\
 & \leq \frac{2\|\nabla\psi\|_{\L^1} \|g_n\circ T_n\|_{\L^p(E)}}{|E|} \frac{\|T_n-\Id\|_{\L^{p^\prime}(E)}}{a_n} \\
 & \qquad + \frac{2\|\nabla\psi\|_{\L^\infty} }{|E|} \frac{\|T_n-\Id\|_{\L^{p^\prime}(E)}}{a_n} \bl \int_E \left|\int_{B(0,1)} g_n(T_n(x-a_nz)) \di z \right|^p \di x \br^{\frac{1}{p}},
\end{align*}
by H\"older's inequality,  where $p^\prime$ satisfies $\frac{1}{p}+\frac{1}{p^\prime}=1$.
Now,
\begin{align*}
\int_E \left|\int_{B(0,1)} g_n(T_n(x-a_nz)) \di z \right|^p \di x & \leq |B(0,1)|^{p-1} \int_E \int_{B(0,1)} \left|g_n(T_n(x-a_nz)) \right|^p \di z \di x \\
 & \leq |B(0,1)|^p \| g_n\circ T_n\|_{\L^p(E)}^p.
\end{align*}
We choose $a_n$ such that, in addition, 
\[ \frac{\|T_n-\Id\|_{\L^{p^\prime}(E)}}{a_n} = O(1)\]
 is satisfied,
and so $\III_n=O(1)$.
The bound on $\IV_n$ follows straightforwardly from
\begin{align*}
\IV_n & = \frac{1}{|E|} \int_E \int_{z\,:\, x-a_nz\in E} |\nabla \psi|(z) |z| \bl g_n(T_n(x)) + g_n(T_n(x-a_n z)) \br \di z \, \di x \\
 & \leq \frac{\|\nabla \psi\|_{\L^\infty}}{|E|} \int_E \int_{B(0,1)} g_n(T_n(x)) + g_n(T_n(x-a_n z)) \di z \di x \\
 & \leq \frac{2\|\nabla \psi\|_{\L^\infty}|B(0,1)| \|g_n\circ T_n\|_{\L^1(E)}}{|E|}.
\end{align*}
Putting the bounds on $\I_n,\II_n,\III_n$ and $\IV_n$ together we can conclude that~\eqref{eq:liminfcompactass1} holds.

To show~\eqref{eq:liminfcompactass2} we write
\begin{align*}
\|& \psi_{a_n} \ast ((c_n\circ T_n) \hat{v}) -(c_n\circ T_n) \hat{v} \|_{\L^1(E)} \\& \leq \int_{E} \left| \int_{\R^d} \psi_{a_n}(x-y) \bl c_n(T_n(y))-c_n(T_n(x))\br \hat{v}(y) \di y \right| \di x \\
 & \qquad \qquad + \int_{E} \left| \int_{\R^d} \psi_{a_n}(x-y) c_n(T_n(x)) \bl \hat{v}(y) - \hat{v}(x) \br \di y \right| \di x \\
 & \leq \int_{E} \int_{\R^d} \psi_{a_n}(x-y) \left| c_n(T_n(y))-c_n(T_n(x)) \right| \hat{v}(y) \di y \di x \\
 & \qquad \qquad + \int_{E} \int_{\R^d} \psi_{a_n}(x-y) \left|c_n(T_n(x))\right| \left| \hat{v}(y) - \hat{v}(x) \right| \di y \di x \\
 & =: \V_n + VI_n.
\end{align*}
By~\eqref{eq:cnDer} we can bound $\V_n$ by
\begin{align*}
\V_n & \leq \frac{1}{|E|} \int_E \int_E \psi_{a_n}(x-y) |T_n(x) - T_n(y)| \bl g_n(T_n(x)) + g_n(T_n(y)) \br \di y \di x \\
 & \leq \frac{2}{|E|} \int_E \int_E \psi_{a_n}(x-y) |T_n(x) - T_n(y)| |g_n(T_n(x))| \di y \di x \\
 & = \frac{2}{|E|} \int_E \int_{z\,:\,x-a_nz\in E} \psi(z) |T_n(x) - T_n(x-a_nz)| |g_n(T_n(x))| \di z \di x \\
 & \leq \frac{2}{|E|} \int_E |T_n(x) - x| |g_n(T_n(x))| \di x \\
 & \qquad \qquad + \frac{2a_n}{|E|} \int_E \int_{z\,:\,x-a_nz\in E} \psi(z) |z| |g_n(T_n(x))| \di z \di x \\
 &  \qquad \qquad + \frac{2}{|E|} \int_{B(0,1)} \int_{\stackrel{w\,:\,w\in E,}{w+a_nz\in E}} \psi(z) |w - T_n(w) |g_n(T_n(w+a_nz))| \di w \di z \\
 & \leq 2\|T_n-\Id\|_{\L^{p^\prime}(\nuv)} \|g_n\|_{\L^p(\nuvn)} + 2a_n\|g_n\|_{\L^1(\nuvn)} \\
  & \qquad \qquad + 2\|T_n-\Id\|_{\L^{p^\prime}(\nuv)}  \int_{B(0,1)} \psi(z)\bl \int_{\stackrel{w\,:\,w\in E,}{w+a_nz\in E}}  |g_n(T_n(w+a_nz)|^p \hat{v}(w) \di w\br^{\frac{1}{p}} \di z.
\end{align*}
Since
\[ 
\int_{B(0,1)} \psi(z)\bl \int_{\stackrel{w\,:\,w\in E,}{w+a_nz\in E}}  |g_n(T_n(w+a_nz)|^p \hat{v}(w) \di w\br^{\frac{1}{p}} \di z\leq \|\psi\|_{\L^\infty} |B(0,1)| \|g_n\|_{\L^p(\nuvn)}
\]
is bounded in $n$,
then $\V_n\to 0$.
The term $\VI_n$ can be bounded as
\begin{align*}
\VI_n & = \int_E |c_n(T_n(x))| \int_{\R^d\setminus E} \psi_{a_n}(x-y) \di y \,\hat{v}(x) \di x \\
 & \leq \| c_n\circ T_n\|_{\L^p(\nuv)} \bl \int_{\dist(x,\partial E)<a_n} \left| \int_{\R^d\setminus E} \psi_{a_n}(x-y) \di y \right|^{p^\prime} \hat{v}(x) \di x \br^{\frac{1}{p^\prime}} \\
 & \leq \|c_n\|_{\L^p(\nuvn)} \bl \frac{|\{x\in E\,:\, \dist(x,\partial E)<a_n\}|}{|E|} \br^{\frac{1}{p^\prime}}
\end{align*}
where, again, $p^\prime$ satisfies $\frac{1}{p}+\frac{1}{p^\prime}=1$. Since $a_n\to 0$ and $\|c_n\|_{\L^p(\nuvn)}$ is bounded, we have $\VI_n\to 0$ by Lemma~\ref{cor:PerBound}. 
Putting the bounds on $\V_n$ and $\VI_n$ together, we conclude that~\eqref{eq:liminfcompactass2} holds.
\end{proof}

\begin{theorem}[Liminf inequality]\label{th:liminfcase1}
	Let $\Omega\subset \R^d$ be an open, bounded set. Let $(v,c^{(1)},c^{(2)})\in CL^p(\Omega)$  and consider positive sequences $\{\epsilon_n\},  \{\mu_{\epsilon_n}\}$, with $\lim_{n\to\infty}\epsilon_n= 0$ and $\lim_{n\to \infty} \mu_{\epsilon_n}\in(0,+\infty]$. Assume that $\{(v_n,c^{(1)}_n,c^{(2)}_n)\}\subset \CL^{p}(\Omega)$ is such that $(v_n,c^{(1)}_n,c^{(2)}_n)\to (v,c^{(1)},c^{(2)})$ in  $\CL^p(\Omega)$. Then, 
	\begin{align*}
	\E_\mu(v,c^{(1)},c^{(2)})\leq \liminf_{n \to \infty} \E_{\mu_{\epsilon_n},\epsilon_n}(v_n,c^{(1)}_n,c^{(2)}_n)
	\end{align*}
	where $\E_{\mu_{\epsilon_n},\epsilon_n}$ and $\E_\mu$ are defined in \eqref{eq:approxenergy} and \eqref{eq:limitenergycase1}, respectively.
\end{theorem}

\begin{proof}
	Since the case $\mu_{\epsilon_n}\to +\infty$ as $n\to \infty$  immediately follows from the case $\mu_{\epsilon_n}\to \mu>0$ as $n\to \infty$, we restrict ourselves to considering $\lim_{n\to \infty}\mu_{\epsilon_n}<+\infty$  in the sequel. 
	Without loss of generality, we can assume that 
	\begin{align*}
	\liminf_{n \to \infty} \E_{\mu_{\epsilon_n},\epsilon_n}(v_n,c^{(1)}_n,c^{(2)}_n)<+\infty,
	\end{align*}
	and by passing to a subsequence (not relabelled) we obtain 
	\begin{align}\label{eq:liminfcond}
	\liminf_{n\to\infty} \E_{\mu_{\epsilon_n},\epsilon_n}(v_n,c^{(1)}_n,c^{(2)}_n)=\lim_{n\to\infty} \E_{\mu_{\epsilon_n},\epsilon_n}(v_n,c^{(1)}_n,c^{(2)}_n)<+\infty.
	\end{align}
	In particular, we can assume, without loss of generality, that $v_n\in W^{1,2}((\Omega,\mathcal{L}^d\lfloor_\Omega);\R),c_n^{(1)}\in W^{1,p}((\Omega,\nuvn);\R^m), c_n^{(2)}\in W^{1,p}((\Omega,\nuvinvn);\R^m)$ for all $n\in\N$.  
	By Theorem~\ref{th:compactness}, the $\CL^p$ convergence of $(v_n,c^{(1)}_n,c^{(2)}_n)\to (v,c^{(1)},c^{(2)})$ and~\eqref{eq:cnl1pliminf}, we have 
	\begin{align*}
	    \liminf_{n\to\infty} \E_{\mu_{\epsilon_n},\epsilon_n}(v_n,c^{(1)}_n,c^{(2)}_n) & = \liminf_{n\to\infty} \Big( \|c_n^{(1)}-u_0\|_{\L^p(\nuvn)}^p + \|c_n^{(2)}-u_0\|_{\L^p(\nuvinvn)}^p \\
	    & \qquad +\mu_{\eps_n} \|c_n^{(1)}\|_{\L^{1,p}(\nuvn)}^p + \mu_{\eps_n} \|c_n^{(2)}\|_{\L^{1,p}(\nuvinvn)}^p + \frac{\nu}{c_W} \E_{\eps_n}^{\GL}(v_n) \Big) \\
	    & \geq \|c^{(1)}-u_0\|_{\L^p(\nuv)}^p + \|c^{(2)}-u_0\|_{\L^p(\nuvinv)}^p + \mu \|c^{(1)}\|_{\L^{1,p}(\nuv)}^p \\
	    & \qquad + \mu \|c^{(2)}\|_{\L^{1,p}(\nuvinv)}^p + \nu \TV(v) \\
	    & = \E_\mu(v,c^{(1)},c^{(2)}),
	\end{align*}
	as required.
\end{proof}

For the limsup inequality we will make use of the following $\L^p$-convergence of translations result.

\begin{proposition}\label{prop:lptranslation}
Let $\Omega\subset \R^d$ be an open, bounded set with Lipschitz boundary, let  $\lambda=\chi_E \mathcal{L}^d\lfloor_{\Omega}$ be the indicator function of some  measurable bounded domain $E\subset \Omega$  with smooth boundary, and let $f\in L^p((\Omega,\lambda);\R^m)$.
Let $\{\lambda_n\}\subset\mathcal{P}(\Omega)$ with Lebesgue densities $\{\rho_n\}\subset L^\infty(\Omega)$. 
Let $S_n\colon \Omega \to \Omega$ be a sequence of transportation maps which pushes forward $\lambda$ to $\lambda_n$, and satisfies $S_n\to \operatorname{Id}$ in $L^p(\Omega;\R^m)$.
Then, 
	\begin{align*}
		\lim_{n\to \infty} \int_\Omega |f(S_n(x))-f(x)|^p \di \lambda(x)=0.
	\end{align*}
\end{proposition}	
\begin{proof}
	Let $\epsilon>0$ be given.  Since $\{\rho_n\}\subset L^\infty(\Omega)$, there exists a constant $C>0$  so that,  for all $n\in \mathbb{N}$, $|\rho_n|\leq C$ $\mathcal{L}^d$-a.e.\ on $\Omega$. Since $f\in L^p((\Omega,\lambda);\R^m)$, we can assume without loss of generality that $f=0$  $\mathcal{L}^d$-a.e.\ on $\Omega\backslash E$. As continuous, compactly supported functions are dense in $L^p(\Omega)$, there exists $g\in C_c(\Omega)$ with $\| f-g\|_{L^{p}(\Omega)}<\tfrac{2\epsilon}{3(1+C^{1/p})}$. 
Further, $$\|f\circ S_n-g \circ S_n\|^p_{L^{p}(\lambda)}= \int_\Omega |f(x) – g(x) |^p  \di \lambda_n(x) = \int_\Omega |f(x) – g(x) |^p \rho_n(x) \di x \leq \left(\frac{2\epsilon}{3(1+C^{1/p})}\right)^pC
 .$$
 For $n\in\N$ sufficiently small, we have  $\|g\circ S_n-g\|_{L^{p}(\Omega)}<\tfrac{\epsilon}{3}$ due to the uniform continuity of $g$. 
 Then, 
	\begin{align*}
		 \bl \int_\Omega |f(S_n(x))-f(x)|^p \di \lambda(x)\br^{\frac{1}{p}} \leq \|f\circ S_n-g \circ S_n\|_{L^{p}(\lambda)}+\|g\circ S_n-g\|_{L^{p}(\lambda)}+\|g-f\|_{L^{p}(\lambda)}<\epsilon,
	\end{align*}
	and this concludes the proof.
\end{proof}	

We now proceed to the limsup inequality.

\begin{theorem}[Limsup inequality]\label{th:limsupcase1}
	Let $\Omega\subset \R^d$ be an open, bounded set with Lipschitz boundary. Let $(v,c^{(1)},c^{(2)})\in CL^p(\Omega)$ with  $\max\{\|c^{(1)}\|_{\L^p},\|c^{(2)}\|_{\L^p}\}<\infty$, and consider positive sequences $\{\epsilon_n\},\{\mu_{\epsilon_n}\}$, with $\lim_{n\to \infty}\epsilon_n= 0$ and $\lim_{n\to \infty} \mu_{\epsilon_n}\in(0,+\infty]$. Then,  there exists a sequence $\{(v_n,c^{(1)}_n,c^{(2)}_n)\}\subset CL^{p}(\Omega)$ such that $(v_n,c^{(1)}_n,c^{(2)}_n)\to (v,c^{(1)}, c^{(2)})$ in  $\CL^p(\Omega)$, and
	\begin{align*}
	\limsup_{n \to \infty} \E_{\mu_{\epsilon_n},\epsilon_n}(v_n,c^{(1)}_n,c^{(2)}_n)\leq \E_\mu(v,c^{(1)},c^{(2)}),
	\end{align*}
	where $\E_{\mu_{\epsilon_n},\epsilon_n}$ and $\E_\mu$ are defined in \eqref{eq:approxenergy} and \eqref{eq:limitenergycase1}, respectively.
\end{theorem}

\begin{proof}
	Without loss of generality, we can assume that $\E_\mu(v,c^{(1)}, c^{(2)})<+\infty$, where $v=\chi_E\in\bv(\Omega;\{0,1\})$ for a measurable set of finite perimeter $E:=\{x\in\Omega\colon v(x)=1\}$, and $c^{(1)}\in \W^{1,p}((\Omega,\nuv);\R^m), c^{(2)}\in \W^{1,p}((\Omega,\nuvinv);\R^m)$. 
	By Theorem~\ref{th:compactness}, there exists  a sequence $\{v_{n}\}\subset \W^{1,2}(\Omega)$ such that $v_{n}\to v$ in $\L^1(\Omega;\R)$ and
	\begin{align*}
	\limsup_{n\to\infty} \int_{\Omega}\left( \epsilon_{n}|\nabla v_{n}|^2+\frac{1}{\epsilon_{n}}W(v_{n})\right) \di x\leq  c_W \tv( v).
	\end{align*}
	We are left to find $c_n^{(1)}\in \W^{1,p}((\Omega,\nuvn);\R^m),c_n^{(2)}\in  \W^{1,p}((\Omega,\nuvinvn);\R^m)$ such that
	\begin{align}
	\limsup_{n\to\infty}\|c_n^{(1)}-u_0\|_{\L^{p}(\nuvn)} &\leq  \|c^{(1)}-u_0\|_{\L^{p}(\nuv)},\label{eq:limsup1}\\ \limsup_{n\to\infty}\|c_n^{(2)}-u_0\|_{\L^{p}(\nuvinvn)} &\leq  \|c^{(2)}-u_0\|_{\L^{p}(\nuvinv)},\\
	\limsup_{n\to\infty}\mu_{\epsilon_n}\|c_n^{(1)}\|^p_{\L^{1,p}(\nuvn)} & \leq \mu	 \|c^{(1)}\|^p_{\L^{1,p}(\nuv)},\label{eq:limsup2}\\ \limsup_{n\to\infty}\mu_{\epsilon_n}\|c_n^{(2)}\|^p_{\L^{1,p}(\nuvinvn)} &\leq \mu	 \|c^{(2)}\|^p_{\L^{1,p}(\nuvinv)},
	\end{align}
	and $(v_n,c^{(1)}_n,c^{(2)}_n)\to (v,c^{(1)}, c^{(2)})$ in $\CL^p(\Omega)$.
	Let $\{T_n^{(1)}\}$ and $\{T_n^{(2)}\}$ be such that ${T_n^{(1)}}_\#\nuv= \nuvn$, ${T_n^{(2)}}_\#\nuvinv= \nuvinvn$, and $\|T_n^{(1)}-\operatorname{Id}\|_{L^p(\nuv)}\to 0$, $\|T_n^{(2)}-\operatorname{Id}\|_{L^p(\nuvinv)}\to 0$, where the existence of $T_n^{(1)},T_n^{(2)}$ is guaranteed by the absolute continuity of $\nuv$ and  $\{\nuvn\}$ converges weakly-$\ast$ to $\nuv$. By Proposition~\ref{prop:convergenceclp}, it suffices to show that
	\begin{align}
	 \lim_{n\to \infty}\|c_n^{(1)}\circ T_n^{(1)} - c^{(1)}\|_{\L^p((\Omega,\nuv);\R^m)}=0,\label{eq:limsup3}\\
	 \lim_{n\to \infty}\|c_n^{(2)}\circ T_n^{(2)}- c^{(2)}\|_{\L^p((\Omega,\nuvinv);\R^m)}=0,
	 \end{align}
	 for $(v_n,c^{(1)}_n,c^{(2)}_n)\to (v,c^{(1)}, c^{(2)})$ in $\CL^p(\Omega)$. 
	
	The proofs for $c^{(2)}$ are analogous to the ones for $c^{(1)}$, so  it  suffices to show the above statements for $c^{(1)}$, i.e., \eqref{eq:limsup1}, \eqref{eq:limsup2}, \eqref{eq:limsup3}. 
	For ease of notation, we drop  the superscript,  write $c$ for $c^{(1)}$, $c_n$ for $c_n^{(1)}$ and $T_n$ for $T_n^{(1)}$, and  assume that $c$ is extended by $0$ on $\R^d\backslash \Omega$.
	
	Let $\psi\in C_c^\infty(\R^d)$ be a standard mollifier (see the proof of Theorem~\ref{th:compactnesssmoothcase1}). 
	We define $c_n:=\psi_{a_n} \ast c\in  C_c^\infty(\Omega;\R^m)$ for any nonnegative, strictly decreasing sequence $\{a_n\}_{n\in\N}\subset \R_+$ with $\lim_{n\to\infty} a_n=0$, which is well-defined due to $\|c\|_{L^p}<\infty$. 
	
	First, we prove \eqref{eq:limsup3}. For this, note that 
	\begin{align*}
	\|c_n\circ T_n-c\|_{\L^{p}(\nuv)}^p
	&= \int_{\Omega} |c_n(T_n(x))-c(x)|^p  \di \nuv(x)\\ 
	&= \int_{\Omega}\left |\int_{\R^d} \psi_{a_n}(T_n(x)-y) (c(y)-c(x))\di y\right |^p  \di \nuv(x)\\ 
	&\leq \int_{\Omega}\int_{\R^d} \psi(z) |c(T_n(x)+a_n z)-c(x)|^p \di z \di \nuv(x)\\
	&= \int_{\R^d} \psi(z) \int_{\Omega}|c(T_n(x)+a_n z)-c(x)|^p  \di \nuv(x)\di z,
	\end{align*}
	where we used the substitution $z:=\tfrac{y-T_n(x)}{a_n}$.
	By the reverse Fatou's Lemma, we obtain
	\begin{align*}
	\limsup_{n\to\infty}\|c_n\circ T_n-c\|_{\L^{p}(\nuv)}^p
	\leq \int_{\R^d} \psi(z) \limsup_{n\to \infty} \int_{\Omega} |c(T_n(x)+a_n z)-c(x)|^p  \di \nuv(x)\di z=0,
	\end{align*}
	where the last equality follows from Proposition \ref{prop:lptranslation}. This yields \eqref{eq:limsup3}.
	
	To show \eqref{eq:limsup1}, note that
	\begin{align*}
	\|c_n-u_0\|_{\L^{p}(\nuvn)} &=  \|c_n\circ T_n-u_0\circ T_n\|_{\L^{p}(\nuv)}\\
	&\leq \|c_n\circ T_n-c\|_{\L^{p}(\nuv)}+\|c- u_0\|_{\L^{p}(\nuv)}+\|u_0- u_0\circ T_n\|_{\L^{p}(\nuv)}.
	\end{align*}
	Hence, \eqref{eq:limsup1} immediately follows from \eqref{eq:limsup3} and Proposition \ref{prop:lptranslation}.

	It remains to prove \eqref{eq:limsup2}. Let  $\{b_n\}_{n\in\N}$ be a sequence  with  $\lim_{n\to\infty} b_n= 0$, whose relation to $\{a_n\}_{n\in\N}$ will be specified below. We introduce the sequence $\{E_{b_n}\}\subset \{x\in E\colon \operatorname{dist}(x,\partial E)>b_n \}$ with smooth boundary, such that $E_{b_n}\to E$ as $n\to \infty$ in the sense that $\chi_{E_{b_n}}\to \chi_E$ in $L^1$ and  $P( E_{b_n})\to P( E)$. 
	For $x,y\in E_{b_n}$, we have
	\begin{align*}
		|c_n(x)-c_n(y)|&=\left|\int_{\R^d} \psi(z) (c(x+a_n z) -c(y+a_n z))\di z\right|\\
		&\leq \int_{\R^d} \psi(z) |c(x+a_n z) -c(y+a_n z)|\di z\\
		&\leq |x-y| \int_{\R^d} \psi(z) \bl g(x+a_n z) +g(y+a_n z)\br\di z\\
		&\leq |x-y| \left(\int_{\R^d} \psi(z) g(x+a_n z)\di z + \int_{\R^d} \psi(z) g(y+a_n z)|\di z\right)\\
		&\leq |x-y| \left((\psi_{a_n}\ast g)(x) +  (\psi_{a_n}\ast g)(y) \right).
	\end{align*}	
	Hence, we obtain
	\begin{align*}
	\|c_n\|_{\L^{1,p}(\nuvn \lfloor_{E_{b_n}})} \leq \|\psi_{a_n}\ast g\|_{\L^{p}(\nuvn \lfloor_{E_{b_n}})}\leq \|\psi_{a_n}\ast g\|_{\L^{p}(\nuvn)}.
	\end{align*}
	Assuming that $g$ is extended by 0 on $\R^d\backslash E$, we have 
	\begin{align*}
		\|\psi_{a_n}\ast g\|_{\L^{p}(\nuvn)}&=\left( \int_\Omega |(\psi_{a_n}\ast g)(x)|^p \di \nuvn(x)\right)^{\frac{1}{p}}
		=\left(\int_\Omega |(\psi_{a_n}\ast g)(T_n(x))|^p \di \nuv(x)\right)^{\frac{1}{p}}\\
		&\leq \left(\int_\Omega |(\psi_{a_n}\ast g)(T_n(x))-g(x)|^p \di \nuv(x)\right)^{\frac{1}{p}}+\| g\|_{\L^{p}(\nuv)}\\
		&\leq \left(\int_{\R^d} \psi(z) \int_\Omega \left| g(T_n(x)+a_n z)-g(x)\right|^p \di \nuv(x)\di z\right)^{\frac{1}{p}}+\| g\|_{\L^{p}(\nuv)},
	\end{align*}
	implying, by Proposition \ref{prop:lptranslation}, that
	\begin{align*}
	\limsup_{n\to\infty} \|c_n\|_{\L^{1,p}(\nuvn \lfloor_{E_{b_n}})} \leq \limsup_{n\to\infty} \|\psi_{a_n}\ast g\|_{\L^{p}(\nuvn)}\leq \| g\|_{\L^{p}(\nuv)}.
	\end{align*}
	Since $\|c\|_{\L^{1,p}(\nuv)}=\inf_g \|g\|_{\L^p(\nuv)}$ by the definition of $\|\cdot\|_{\L^{1,p}(\nuv)}$, this yields
		\begin{align*}
	\limsup_{n\to\infty} \|c_n\|_{\L^{1,p}(\nuvn \lfloor_{E_{b_n}})}\leq \| c\|_{\L^{1,p}(\nuv)}.
	\end{align*}
	We denote the complement of $E_{b_n}$ in $\R^d$ by $E_{b_n}^c$ and, since $c_n\in C^\infty$, we have 
	\begin{align*}
		\|c_n\|_{\L^{1,p}(\nuvn)}= \|c_n\|_{\L^{1,p}(\nuvn \lfloor_{E_{b_n}})} + \|c_n\|_{\L^{1,p}(\nuvn \lfloor_{E_{b_n}^c})}.
	\end{align*}
	It remains to show that
	\begin{align*}
	\limsup_{n\to\infty}\|c_n\|_{\L^{1,p}(\nuvn \lfloor_{E_{b_n}^c})}= 0.
	\end{align*}
	We have 
	\begin{align*}
		\|c_n\|_{\L^{1,p}(\nuvn \lfloor_{E_{b_n}^c})}^p & = \|\nabla c_n\|_{\L^{p}(\nuvn \lfloor_{E_{b_n}^c})}^p = \int_{\Omega\backslash E_{b_n}} |\nabla c_n|^p \di \nuvn(x) \\ 
		& = \int_{\Omega\backslash E_{b_n}} |(\nabla \psi_{a_n}\ast c)(T_n(x))|^p \di \nuv(x) \\ 
		& = \int_{\Omega\backslash E_{b_n}} \left|\int_{\R^d} \frac{1}{a_n^{d+1}}\nabla \psi\left( \frac{T_n(x)-y}{a_n}\right) c(y)\di y\right|^p \di \nuv(x) \\ 
		& = \int_{\Omega\backslash E_{b_n}} \left|\int_{\R^d} \frac{1}{a_n}\nabla \psi(z) c(T_n(x)-a_n z)\di z\right|^p \di \nuv(x) \\
		& \leq \frac{1}{a_n^p} \|c\|_{\L^p(\Omega)}^p\|\nabla \psi\|_{\L^\infty}^p \int_{\Omega\backslash E_{b_n}} \di \nuv(x).
	\end{align*}
	Suppose that
	\begin{align*}
		a_n:=\left(\int_{\Omega\backslash E_{b_n}} \di \nuv(x)\right)^\frac{1}{2p},\quad n\in \N,
	\end{align*}
	so that $\lim_{n\to\infty} a_n=0$ as required above. Then, 
		\begin{align*}
	\limsup_{n\to\infty}\|c_n\|_{\L^{1,p}(\nuvn \lfloor_{E_{b_n}^c})}= 0,
	\end{align*}
	which yields
	\begin{align*}
	\limsup_{n\to\infty} \|c_n\|_{\L^{1,p}(\nuvn)}\leq \| c\|_{\L^{1,p}(\nuv)}.
	\end{align*}
	
	If $\lim_{n\to \infty} \mu_{\epsilon_n}=\mu>0$, then we have 
	\begin{align*}
	\limsup_{n\to\infty}\mu_{\epsilon_n}\|c_n\|_{\L^{1,p}(\nuvn)}\leq \mu	 \|c\|_{\L^{1,p}(\nuv)},
	\end{align*}
	which concludes the limsup inequality.
	
	For $\lim_{n\to\infty} \mu_{\epsilon_n}=+\infty$, 
	$c$ is constant $\mathcal{L}^d$-a.e.\ $x\in E$.
	This implies that
	\[ \limsup_{n\to\infty}	\|c_n\|^p_{\L^{1,p}(\nuvn)}=\|c\|^p_{\L^{1,p}(\nuv)} =0, \]
	and hence the limsup inequality also holds for $\lim_{n\to\infty} \mu_{\epsilon_n}=+\infty$.
\end{proof}

The $\Gamma$-convergence result in Theorem \ref{th:mainresult} follows from the liminf inequality in Theorem \ref{th:liminfcase1} and the limsup inequality in Theorem \ref{th:limsupcase1}. Note that the property  $\max\{\|c^{(1)}\|_{\L^p},\|c^{(2)}\|_{\L^p}\}<\infty$ in Theorem \ref{th:limsupcase1} is used to simplify the notation as for any $(v,c^{(1)}, c^{(2)})\in \CL^p(\Omega)$ we have $c^{(1)}\in \L^p((\Omega,\nuv);\R^m), c^{(2)}\in \L^p((\Omega,\nuvinv);\R^m)$ and hence we can assume without loss of generality that $\max\{\|c^{(1)}\|_{\L^p},\|c^{(2)}\|_{\L^p}\}<\infty$ holds.

Due to the compactness property  in Theorem \ref{th:compactnesssmoothcase1} with  regularity assumptions on $E$ and the $\Gamma$-convergence of the energy functionals,   we can conclude the convergence of minimizers $(v_n,c_n^{(1)},c_n^{(2)})$, see Corollary \ref{cor:convergenceminimizer}, once we have shown that $\sup_{n\in\N}\max\{\|c_n^{(1)}\|_{\L^\infty},\|c_n^{(2)}\|_{\L^\infty}\}<\infty$.

\begin{proof}[Proof of Corollary~\ref{cor:convergenceminimizer}]
To show that $\sup_{n\in\N}\max\{\|c_n^{(1)}\|_{\L^\infty},\|c_n^{(2)}\|_{\L^\infty}\}<\infty$, 
 we suppose  $m=1$ for simplicity, i.e., $u_0:\Omega\to\R$. One can proceed in a similar way for $m>1$.
Let $M:=\|u_0\|_{\L^\infty}$, and assume that $(v_n,c_n^{(1)},c_n^{(2)})$ is a minimizer of $\E_{\mu_{\eps_n},\eps_n}$.
For a contradiction, we suppose that there exists $i\in \{1,2\}$ and $n\in \N$ such that $\|c_n^{(i)}\|_{\L^\infty}>M+1$.
We define
\begin{align*}
\tilde{c}^{(i)}_n(x) := \left\{ \begin{array}{ll} M, & \text{if } c_n^{(i)}(x) > M, \\ c_n^{(i)}(x), & \text{if } c_n^{(i)}(x)\in [-M,M], \\ -M, & \text{if } c_n^{(i)}(x)<-M. \end{array} \right.
\end{align*}
Clearly $\|\tilde{c}_n\|_{\L^\infty}\leq M$.
Moreover,
\begin{align*}
\|\tilde{c}_n^{(i)} - u_0\|_{\L^p(\nuvn)}^p & = \int_{|c_n^{(i)}(x)|\leq M} \left| c_n^{(i)} - u_0(x)\right|^p \di \nuvn(x) \\
 & \qquad + \int_{c_n^{(i)}(x)\in (-M-1,- M)} \underbrace{\left| -M - u_0(x)\right|^p}_{\leq |c_n^{(i)}(x) - u_0(x)|^p} \di \nuvn(x) \\
 & \qquad + \int_{c_n^{(i)}(x)\in (M,M+1)} \underbrace{\left| M - u_0(x)\right|^p}_{\leq |c_n^{(i)}(x) - u_0(x)|^p} \di \nuvn(x) \\
 & \qquad + \int_{c_n^{(i)}(x)\leq -M-1} \underbrace{\left| -M - u_0(x)\right|^p}_{\leq |c_n^{(i)}(x) - u_0(x)|^p - 1} \di \nuvn(x) \\
 & \qquad + \int_{c_n^{(i)}(x)\in \geq M+1} \underbrace{\left| M - u_0(x)\right|^p}_{\leq |c_n^{(i)}(x) - u_0(x)|^p-1} \di \nuvn(x) \\
 & \leq \| c_n^{(i)} - u_0\|_{\L^p(\nuvn)}^p - \underbrace{\left| \left\{ x\,:\, |c_n^{(i)}(x)|>M+1 \right\} \right|}_{>0}.
\end{align*}
One can easily check that, for all $x,y$,
\[ \left| \tilde{c}_n^{(i)}(x) - \tilde{c}_n^{(i)}(y) \right| \leq \left| c_n^{(i)}(x) - c_n^{(i)}(y) \right|, \]
and therefore $\|\tilde{c}_n^{(i)}\|_{\L^{1,p}(\nuvn)} \leq \|c_n^{(i)}\|_{\L^{1,p}(\nuvn)}$.
We have shown that
\[ \E_{\mu_{\eps_n},\eps_n}(v_n,\tilde{c}_n^{(1)},\tilde{c}_n^{(2)}) < \E_{\mu_{\eps_n},\eps_n}(v_n,c_n^{(1)},c_n^{(2)}), \] 
which contradicts the assumption that $(v_n,c_n^{(1)},c_n^{(2)})$ is a minimizer.
Hence, $\|c_n^{(i)}\|_{\L^\infty}<M+1$ for all $i=1,2$, and $n\in \N$.
\end{proof}

\section*{Acknowledgements}
The authors thank Francesco Maggi for his advice and references on  isoperimetric inequalities. I.\ Fonseca acknowledges the Center for Nonlinear Analysis (CNA) where part of this work was carried out. Her research was partial funded under grants NSF DMS No.\ 1411646, No.\ 1906238 and No.\ 2205627. L.~M.~Kreusser, C.-B.\ Schönlieb and M.~Thorpe would like to thank the Isaac Newton Institute for Mathematical Sciences for support and hospitality during the programme \emph{Mathematics of Deep Learning} when work on this paper was undertaken (EPSRC grant number EP/R014604/1). L.~M.~Kreusser, C.-B.\ Schönlieb and M.~Thorpe acknowledge support from the European Union Horizon 2020 research and innovation programmes under the Marie Sk\l odowska-Curie grant agreement No. 777826 (NoMADS). L.~M.~Kreusser also acknowledges support the EPSRC grant EP/L016516/1, the German National Academic Foundation (Studienstiftung des Deutschen Volkes), the Cantab Capital Institute for the Mathematics of Information and Magdalene College, Cambridge (Nevile Research Fellowship). C.-B.\ Schönlieb acknowledges support from the Philip Leverhulme Prize, the Royal Society Wolfson Fellowship, the EPSRC advanced career fellowship EP/V029428/1, EPSRC grants EP/S026045/1 and EP/T003553/1, EP/N014588/1, EP/T017961/1, the Wellcome Innovator Award RG98755, the Cantab Capital Institute for the Mathematics of Information and the Alan Turing Institute.
M.~Thorpe also holds a Turing Fellowship at the Alan Turing Institute.

\bibliographystyle{plain}
\bibliography{references}

\appendix

\section{Enlarged Boundaries for Sets of Finite Perimeter}

For completeness, we include a bound on the volume
$$ \left| \left\{ x\in \R^d \,:\, \dist(x,\partial E)\leq a \right\}\right|, $$
where  $E$ denotes a set with finite perimeter which is used in the compactness result in Theorem~\ref{th:compactnesssmoothcase1}.

\begin{theorem}\cite[Theorem 2.106]{ambrosio2000functions}\label{th:ambrosio}
If $Z$ is a compact, countably $H^k$-rectifiable set in $\R^d$ and if there are $\kappa>0$ and $r_0>0$ such that
$H^k(Z\cap B_r(x))\ge \kappa r^k $         for every $x \in Z$ and every $ r<r_0$,
then $Z$ is $k$-Minkowski regular, i.e.\ there exists a constant $o>0$ such that
\begin{align*}
|\{x\in \R^d\colon \dist(x,Z)< a\}|  =  om_{d-k} a^{d-k}  H^k(Z)  + o(a^{d-k})  \quad \text{as} \quad   a\to 0,
\end{align*}
where $m_{d-k}$ denotes the $d-m$ dimensional sphere.
\end{theorem}

Applying Theorem \ref{th:ambrosio} to our setting yields an estimate for $\left| \left\{ x\in \R^d \,:\, \dist(x,\partial E)\leq a \right\}\right|$:
\begin{corollary}\label{cor:PerBound}
  Let $E$ be a bounded set of finite perimeter in $\R^{d}$. Assume  that the topological boundary $\partial E$ is the closure of the reduced boundary $\partial^* E$.  
Assume that for some $\kappa>0$ and some $r_0>0$ we have $P(E;B_r(x)) \geq \kappa r^{d-1}$ for every $x \in \partial^* E$.
Then $$|\{x \in \R^d\colon \dist(x,\partial E)< a\}|  =  2 a  P(E)  + o(a)     \quad \text{as} \quad     a\to 0.$$
\end{corollary}

\begin{proof}
The assumptions on $E$ imply that $E$ is compact and countably $H^{d-1}$-rectifiable.
 Since $\partial E$ is the closure of $\partial^* E$, the fact that $P(E;B_r(x))=H^{d-1}(B_r(x)\cap\partial^* E)$ has lower density estimates implies by continuity that $Z=\partial E$ has lower density estimates, and then one applies Theorem \ref{th:ambrosio}.
\end{proof}

Note that the assumption on the topological boundary in Corollary \ref{cor:PerBound} holds for free up to a modification on a Lebesgue null set, see \cite[Proposition 12.20]{maggi2012sets}. Hence, when applying Corollary \ref{cor:PerBound} to a bounded set $E$ with finite perimeter in the proof of Theorem \ref{th:compactnesssmoothcase1}, it is sufficient to assume that for some $\kappa>0$ and some $r_0>0$ we have $P(E;B_r(x)) \geq \kappa r^{d-1}$ for every $x \in \partial^* E$.

\end{document}